\documentclass[11pt]{article}

\usepackage{latexsym,epsfig}
\usepackage{amssymb,amsmath,amsfonts}
\usepackage{exscale}
\usepackage{ifthen}
\usepackage{longtable}
\usepackage{subfigure}
\usepackage{caption}
\usepackage{cite}
\usepackage{lscape}
\usepackage{hyperref}

\begin{document}
\bibliographystyle{gabialpha}
\protect\pagenumbering{arabic}
\setcounter{page}{1}

\newcommand{\Zt}{\rm}

\newcommand{\ba}{\begin{array}}
\newcommand{\ea}{\end{array}}
\newcommand{\pot}{{\cal P}}
\newcommand{\curv}{\cal C}
\newcommand{\ddt} {\mbox{$\frac{\partial  }{\partial t}$}}
\newcommand{\hl}{\sf}
\newcommand{\hd}{\sf}

\newcommand{\Ad}{\mbox{\rm Ad}}
\newcommand{\Adsm}{\mbox{{\rm \scriptsize Ad}}}
\newcommand{\ad}{\mbox{\rm ad}}
\newcommand{\adsm}{\mbox{{\rm \scriptsize ad}}}
\newcommand{\diag}{\mbox{\rm Diag}}
\newcommand{\sect}{\mbox{\rm sec}}
\newcommand{\id}{\mbox{\rm id}}
\newcommand{\idsm}{\mbox{{\rm \scriptsize id}}}
\newcommand{\eps}{\varepsilon}

\newcommand{\aL}{\mathfrak{a}}
\newcommand{\bL}{\mathfrak{b}}
\newcommand{\mL}{\mathfrak{m}}
\newcommand{\kL}{\mathfrak{k}}
\newcommand{\gL}{\mathfrak{g}}
\newcommand{\nL}{\mathfrak{n}}
\newcommand{\hL}{\mathfrak{h}}
\newcommand{\pL}{\mathfrak{p}}
\newcommand{\uL}{\mathfrak{u}}
\newcommand{\lL}{\mathfrak{l}}

\newcommand{\kG}{{\tt k}}
\newcommand{\nG}{{\tt n}}

\newcommand{\Cart}{$G=K e^{\overline{\aL^+}} K$}
\newcommand{\Area}{\mbox{Area}}
\newcommand{\Hd}{\mbox{\rm Hd}}
\newcommand{\Hdim}{\mbox{\rm dim}_{\mbox{\rm \scriptsize Hd}}}
\newcommand{\Tr}{\mbox{\rm Tr}}
\newcommand{\bs}{{\cal B}}
\newcommand{\bv}{{\rm B}}

\newcommand{\nc}{{\cal N}}
\newcommand{\MM}{{\cal M}}
\newcommand{\Ch}{{\cal C}}
\newcommand{\clCh}{\overline{\cal C}}
\newcommand{\Sh}{\mbox{Sh}}
\newcommand{\smSh}{\mbox{{\rm \scriptsize Sh}}}
\newcommand{\Cnt}{\mbox{\rm C}}
\newcommand{\preim}{(\pi^F)^{-1}}

\newcommand{\NN}{\mathbb{N}} \newcommand{\ZZ}{\mathbb{Z}}
\newcommand{\QQ}{\mathbb{Q}} \newcommand{\RR}{\mathbb{R}}
\newcommand{\KK}{\mathbb{K}} \newcommand{\FF}{\mathbb{F}}
\newcommand{\CC}{\mathbb{C}} \newcommand{\EE}{\mathbb{E}}
\newcommand{\XX}{X}
\newcommand{\HH}{I\hspace{-2mm}H}
\newcommand{\norm}{\Vert\hspace{-0.35mm}|}
\newcommand{\Sph}{\mathbb{S}}
\newcommand{\ganz}{\overline{\XX}}
\newcommand{\rand}{\partial\XX}
\newcommand{\prodrand}{\partial\XX_1\times\partial\XX_2} 
\newcommand{\regrand}{\partial\XX^{reg}}
\newcommand{\singrand}{\partial\XX^{sing}}
\newcommand{\Frand}{\partial^F\XX}
\newcommand{\Lim}{L_\Gamma}          
\newcommand{\Flim}{F_\Gamma}
\newcommand{\reglim}{L_\Gamma^{reg}}
\newcommand{\radlim}{L_\Gamma^{rad}}
\newcommand{\raylim}{L_\Gamma^{ray}}
\newcommand{\horinf}{\mbox{Vis}^{\infty}}
\newcommand{\horF}{\mbox{Vis}^F}
\newcommand{\Sml}{\mbox{Small}}
\newcommand{\SmlF}{\mbox{Small}^F}

\newcommand{\ifl}{\qquad\Longleftrightarrow\qquad}
\newcommand{\at}{\!\cdot\!}
\newcommand{\ging}{\gamma\in\Gamma}
\newcommand{\xo}{{o}}
\newcommand{\gamo}{{\gamma\xo}}
\newcommand{\gam}{\gamma}
\newcommand{\gax}{h}
\newcommand{\gxi}{{G\!\cdot\!\xi}}
\newcommand{\bd}{$(b,\theta)$-densit}
\newcommand{\bt}{$(b,\theta)$-densit}
\newcommand{\cd}{$(\alpha,\Gamma\at\xi)$-density}
\newcommand{\be}{\begin{eqnarray*}}
\newcommand{\ee}{\end{eqnarray*}}

\newcommand{\e}{\mathrm{e}}
\newcommand{\x}{p}
\newcommand{\y}{q}

\newcommand{\an}{\ \mbox{and}\ }
\newcommand{\as}{\ \mbox{as}\ }
\newcommand{\diam}{\mbox{diam}}
\newcommand{\is}{\mbox{Is}}
\newcommand{\Ax}{\mbox{Ax}}
\newcommand{\Fix}{\mbox{Fix}}
\newcommand{\Par}{F}
\newcommand{\Min}{\mbox{Fix}}
\newcommand{\rel}{\mbox{Rel}_\Gamma}
\newcommand{\vol}{\mbox{vol}}
\newcommand{\Td}{\mbox{Td}}
\newcommand{\piF}{\pi^B}
\newcommand{\piKM}{\pi^I}

\newcommand{\for}{\ \mbox{for}\ }
\newcommand{\pr}{\pi}
\newcommand{\sh}{\mbox{sh}}
\newcommand{\shi}{\mbox{sh}^{\infty}}
\newcommand{\rank}{\mbox{rank}}
\newcommand{\supp}{\mbox{supp}}
\newcommand{\mass}{\mbox{mass}}
\newcommand{\kernel}{\mbox{kernel}}
\newcommand{\st}{\mbox{such}\ \mbox{that}\ }
\newcommand{\Stab}{\mbox{Stab}}
\newcommand{\Root}{\Sigma}
\newcommand{\Cone}{\mbox{C}}
\newcommand{\wrt}{\mbox{with}\ \mbox{respect}\ \mbox{to}\ }
\newcommand{\where}{\ \mbox{where}\ }

\newcommand{\thet}{\widehat H}

\newcommand{\con}{{\sc Consequence}\newline}
\newcommand{\rem}{{\sc Remark}\newline}
\newcommand{\prf}{{\sl Proof.\  }}
\newcommand{\qed}{$\hfill\Box$}

\newenvironment{rmk} {\newline{\sc Remark.\ }}{}  
\newenvironment{rmke} {{\sc Remark.\ }}{}  
\newenvironment{rmks} {{\sc Remarks.\ }}{}  
\newenvironment{nt} {{\sc Notation}}{}  

\newtheorem{satz}{\bf Theorem}

\newtheorem{df}{\sc Definition}[section]
\newtheorem{cor}[df]{\sc Corollary}
\newtheorem{thr}[df]{\bf Theorem}
\newtheorem{lem}[df]{\sc Lemma}
\newtheorem{prp}[df]{\sc Proposition}
\newtheorem{ex}{\sc Example}
\newenvironment{pros}{{\sc Properties:}}


\title{{\sc Asymptotic geometry in higher products of rank one Hadamard spaces}}
\author{\sc Gabriele Link}
\date{\today}
\maketitle
\begin{abstract} Given a product $\XX$ 
of locally compact rank one Hadamard spaces, we study asymptotic properties of certain discrete isometry groups $\Gamma$ of $\XX$. 
First we give a detailed description of the structure of the geometric limit set and relate it to the limit cone; moreover, we show that the action  of $\Gamma$ on a quotient of the regular geometric boundary of $\XX$ is minimal and proximal. This is completely analogous to the case of Zariski dense discrete subgroups of semi-simple Lie groups acting on the associated symmetric space (compare \cite{MR1437472}).  In the second part of the paper we study the distribution of $\Gamma$-orbit points in $\XX$: As a generalization of the critical exponent $\delta(\Gamma)$ of $\Gamma$ we consider for any $\theta\in\RR_{\ge 0}^r$, $\Vert\theta\Vert=1$, the exponential growth rate  $\delta_\theta(\Gamma)$  of the number of orbit points in $\XX$  with prescribed {``}slope" $\theta$.  
In analogy to Quint's result in  \cite{MR1933790} we show that  the  homogeneous extension $\Psi_\Gamma$ to $\RR_{\ge 0}^r$ of  $\delta_\theta(\Gamma)$ as a function of $\theta$ is upper semi-continuous, concave and strictly positive in the relative interior of the limit cone. This shows in particular that there exists a unique slope $\theta^*$ for which $\delta_{\theta^*}(\Gamma)$ is maximal and equal to the critical exponent 
of $\Gamma$.  

We notice that an interesting class of product spaces as above comes from the second alternative in 
the Rank Rigidity Theorem (\cite[Theorem A]{MR2827012})  for 
CAT$(0)$-cube complexes: 
Given a finite-dimensional CAT$(0)$-cube complex $\XX$ and a group $\Gamma$ of automorphisms without fixed point in the geometric compactification of $\XX$, then either $\Gamma$ contains a rank one isometry or there exists a convex $\Gamma$-invariant subcomplex of $\XX$ which is a product of two unbounded cube subcomplexes; in the latter case one inductively gets a convex $\Gamma$-invariant subcomplex of $\XX$ which can be decomposed into a finite product of rank one Hadamard spaces. 
So our results  imply in particular that 
classical properties of discrete subgroups of higher rank Lie groups as in  \cite{MR1437472} and \cite{MR1933790} also hold for certain discrete isometry groups of reducible CAT$(0)$-cube complexes.  

\end{abstract}

\vspace{0.2cm}

\section{Introduction}
After the publication of the article {``}Asymptotic geometry in products of Hadamard spaces with rank one isometries" (\cite{MR2629900})  I was asked several times whether the results naturally extend to the setting of more than two factors. Unfortunately this is not the case since the methods of proof used there rely heavily on the possibility to control the position of pairs of points in a quotient of the regular geometric boundary; when more factors are present, the 
set of pairs of points which are in an uncontrollable position becomes much larger -- in the case of symmetric spaces this phenomenon is reflected in the presence of more and higher dimensional {``}small" Bruhat cells when the rank gets bigger. So one goal of the present article was to give a generalization of the results in the aforementioned paper to products with more than two factors. Apart from that, the article contains a variety of results which were not yet known in the case of two factors. Among these I only want to mention  here the construction of freely generated discrete subgroups (Proposition~\ref{constrfreegrps}) with limit cone contained in a prescribed set and the positivity of the critical exponent (Theorem~\ref{expgrowthpositive}).

To be more precise, 
we let $(\XX,d)$ be a product of $r$ locally compact Hadamard spaces $(\XX_i,d_i)$ endowed with the $\ell^2$-metric, which makes $\XX$ itself a locally compact 
Hadamard space, i.e. a locally compact complete simply connected metric space of non-positive Alexandrov curvature. 
It is well-known that every  locally compact Hadamard space $\XX$ can be 
compactified by adding its  geometric boundary $\rand$ 
endowed with the cone topology (see \cite[chapter II]{MR1377265});
if $\XX$ is a product, then 
the {\hl regular geometric boundary} $\regrand$ of $\XX$ -- which consists of the set of equivalence classes of geodesic rays which do not project to a point in one of the factors -- is a dense open subset of $\rand$ homeomorphic to the Cartesian product of the geometric boundaries $\rand_i$ of the factors $\XX_i$ (which we call the {\hl Furstenberg boundary} $\Frand$ of $\XX$) times a factor 
$E^+=\{\theta\in\RR_{> 0}^r: \Vert\theta\Vert=1\}$.


We next let $\Gamma<\is(\XX)$  be a group acting properly discontinuously by isometries on $\XX$; passing to a subgroup of finite index if necessary we may further assume that $\Gamma$ preserves the product decomposition (\cite[Corollary~1.3]{MR2399098}).  The {\hl geometric limit set} $\Lim\subset\rand$ of $\Gamma$ is defined as the set of accumulation points of a $\Gamma$-orbit in $\XX$. 
Unlike in the case of  CAT$(-1)$-spaces this geometric limit set is in general not  a minimal set for the action of $\Gamma$ on the geometric boundary $\rand$ of $\XX$: This is due to the fact that isometries preserving the product decomposition of $\XX$ cannot change the {\hl slope} (i.e. the projection  to $E^+$) of regular boundary points.

In this note we further restrict our attention to discrete groups $\Gamma$ as above 
which contain an element projecting to a rank one isometry in each factor, i.e. $\Gamma$ contains an element $h$ 
\st all its projections $h_i$ to $\is(\XX_i)$ possess an 
invariant geodesic 
which does not bound a flat half-plane in the corresponding factor $\XX_i$.
This requires in particular that all factors are rank one, i.e. possess a geodesic line which does not bound a flat half-plane. 
For products of rank one  Hadamard spaces the presence of such a {\hl regular axial} isometry is guaranteed in many interesting cases: If for example $\XX$ is a finite-dimensional CAT$(0)$-cube complex for which every irreducible factor is 
non-Euclidean,  unbounded and locally compact with a cocompact and essential action of its automorphism group, then by Theorem C in \cite{MR2827012} every (possibly non-uniform) lattice $\Gamma<\is(\XX)$ contains a regular axial isometry.

We will moreover need a second regular axial isometry $g\in \Gamma$ 
\st all projections to $\is(\XX_i)$ of $g$ and $h$ are independent. 
 This  condition is clearly satisfied when $\Gamma$ contains a (not necessarily regular axial) element $\gamma$ 
 \st all projections 
$\gamma_i$ to $\is(\XX_i)$ map the two fixed points of $h_i$ to their complement in $\rand_i$; then 
$g= \gamma h\gamma^{-1}\in\Gamma$ is the desired second regular axial isometry.
Another important class of examples satisfying this assumption in the case of only two factors are Kac-Moody groups  $\Gamma$ over a finite field: They act by isometries on a product, 
the CAT$(0)$-realization of the associated twin building ${\cal B}_+\times{\cal B}_-$, and there  exists an element $h=(h_1,h_2)\in\Gamma$ projecting to a rank one element in each factor by  Remark 5.4 and the proof of Corollary~1.3 in \cite{MR2585575}. Moreover, the action of the Weyl group produces many regular axial isometries $g=(g_1,g_2)\in\Gamma$ with $g_1$ independent from $h_1$ and $g_2$ independent from $h_2$. 
Notice that if the order of the ground field is sufficiently large, then $\Gamma$
is an irreducible lattice (see e.g. \cite{MR1715140} and \cite{MR2485882}). 

Apart from these examples possible factors of $\XX$ include locally compact Hadamard spaces of strictly negative Alexandrov curvature such as locally finite trees or manifolds with sectional curvature bounded from above by a negative constant as e.g. rank one symmetric spaces of the non-compact type.  
In this special case every non-elliptic and non-parabolic isometry in one of the factors is already a rank one element.  Prominent examples here which are already covered by the results of Y.~Benoist and J.-F.~Quint are Hilbert modular groups acting as irreducible lattices on a product of hyperbolic planes, and graphs of convex cocompact groups of rank one symmetric spaces (see also \cite{MR1230298}). More generally, given a set of locally compact rank one Hadamard spaces $\XX_1$, $\XX_2 ,\ldots,\XX_r$, and faithful representions $\rho_i:G\to\is(\XX_i)$  of a group $G$ acting properly discontinuously by isometries on $\XX_i$, $i\in\{1,2,\ldots,r\}$,  the {\hl graph group} $\Gamma_\rho$ associated to $\rho=(\rho_1,\rho_2,\ldots,\rho_r)$ is defined by
\[\Gamma_\rho=\{(\rho_1(g),\rho_2(g),\ldots,\rho_r(g)):g\in G\};\]
it clearly satisfies our assumption if $G$ contains two elements $g$ and $h$ \st $\rho_i(g)$ and $\rho_i(h)$ are independent rank one isometries of $\XX_i$ for all $i\in\{1,2,\ldots,r\}$. Notice that thanks to general arguments given by F.~Dal'bo in \cite[Proposition~3.4]{MR1604251} 
any discrete group $\Gamma<\is(\XX_1)\times\is(\XX_2)\times\cdots\times \is(\XX_r)$  which acts freely 
and  
satisfies 
$\Lim\subset\regrand$\break is a graph group; the converse clearly does not hold.

We will now state our main results. The first one 
is a strengthening of Theorem~A in \cite{MR2629900}:\\[2mm]
{\bf Theorem A}$\quad$ {\sl The Furstenberg limit set is a boundary limit set for the action of $\,\Gamma$ on the Furstenberg boundary $\Frand$. }\\[3mm]
This implies in particular that the Furstenberg limit set is minimal, i.e. the smallest non-empty, $\Gamma$-invariant closed subset of $\Frand$. Notice that in the recent article \cite{1105.1675}  A.~Nevo and M.~Sageev proved an analogous statement  for the Poisson boundary $B(\XX)$ of a proper cocompact $G$-action on an essential, strictly non-Euclidean CAT$(0)$-cube complex $\XX$. More precisely, their Theorem~5.8 states that the closure $B(\XX)$ of the set of non-terminating ultrafilters  in the Roller boundary of $\XX$ 
is a boundary limit set for the action of $G$ on the collection 
of all ultrafilters. 

In our setting we have -- as in the case of symmetric spaces or Bruhat-Tits buildings of higher rank -- the following structure theorem:\\[2mm]
{\bf Theorem B}$\quad$ {\sl The regular geometric limit set splits as a product $\Flim\times P_\Gamma^{reg}$, where $P_\Gamma^{reg}\subset E^+$ denotes the set of  slopes of  regular limit points of $\,\Gamma$. }\\[3mm]
We recall that 
$\Gamma<\is(\XX_1)\times\is(\XX_2)\times\cdots\times\is(\XX_r)$
is a discrete group 
containing two  regular axial isometries $h=(h_1,h_2,\ldots,h_r)$ and $g=(g_1,g_2,\ldots,g_r)$ which project to independent rank one elements in each factor.  For a rank one isometry $\gamma_i$ of one of the factors $\XX_i$ we denote $\gamma_i^+$ its attractive, $\gamma_i^-$ its repulsive fixed point in $\rand_i$, and $l_i(\gamma_i)>0$ 
its translation length, i.e. the minimum of the displacement function $d_i(\gamma_i):\XX_i\to\RR$, $x_i\mapsto d_i(x_i,\gamma_i x_i)$. So if $\gamma=(\gamma_1,\gamma_2,\ldots,\gamma_r)$ is regular axial we canonically get two fixed points in the Furstenberg boundary $\gamma^+:=(\gamma_1^+,\gamma_2^+,\ldots,\gamma_r^+)$,     
$\gamma^-:=(\gamma_1^-,\gamma_2^-,\ldots,\gamma_r^-)$, 
and a {\hl translation vector}  
$L(\gamma):=(l_1(\gamma_1),l_2(\gamma_2),\ldots,l_r(\gamma_r))\in\RR_{>0}^r$.
With this notation we can state the following \\[2mm] 
{\bf Theorem C}$\quad$ {\sl The set of pairs of fixed points $(\gamma^+,\gamma^-)\in\Frand\times \Frand$ of regular axial isometries $\gamma\in\Gamma$ is dense in $\big(F_\Gamma
\times F_\Gamma\big)\setminus \Delta$, where $\Delta$ denotes the set of pairs of points in $F_\Gamma$ with a common projection to some $\rand_i$.  }\\[3mm] 
We mention that this result can be viewed as a  strong topological version of the double ergodicity property of Poisson boundaries due to Burger-Monod (\cite{MR1911660}) and Kaimanovich (\cite{MR2006560}). 

We next define the {\hl limit cone} $\ell_\Gamma$ of $\Gamma$ as the closure in $\RR_{\ge 0}^r$ of the set of half-lines  spanned by 
all translation vectors $L(\gamma)$ 
with $\gamma\in\Gamma$ regular axial. 
This cone is related to the 
set of slopes of all (regular and singular) limit points as follows:\\[2mm]  
{\bf Theorem D}$\quad$ {\sl The set $P_\Gamma\subset E$ of  slopes of  all limit points of $\,\Gamma$ 
is equal to $\ell_\Gamma\cap E$. 
Moreover, the limit cone $\ell_\Gamma$ 
is convex.}\\[3mm]
For the formulation of the last main result of the paper we fix a point $x=(x_1,x_2,\ldots, x_r)$ in $\XX$, 
$\theta=(\theta_1,\theta_2,\ldots,\theta_r)\in E:= \overline{E^+}\subset\RR_{\ge 0}^r$, $\eps>0$, $n\gg 1$  and consider the cardinality $N_\theta^\eps(n)$ of  the set
\begin{align*} \big\{ \gamma=(\gamma_1,\gamma_2,\ldots,\gamma_r)\in\Gamma: &\  0<d(x,\gamma x)<n,\\ 
&\  \Big| \frac{d_i(x_i,\gamma_i x_i)}{d(x,\gamma x)}-\theta_i\Big|<\eps\ \text{ for all }\  1\le i\le r\}\big\}.
\end{align*} 
This number counts all orbit points $\gamma x$ 
of distance less than $n$ to the point $x$ which in addition are {``}close" to a geodesic ray in the class of a boundary point with slope $\theta$. 
We further consider 
$$\delta_\theta^\eps:=\limsup_{n\to\infty}\frac{\log N_\theta^\eps(n)}{n} \quad\text{ and finally }\quad \delta_\theta(\Gamma):=\lim_{\eps\to 0} \delta_\theta^\eps.$$
$\delta_\theta(\Gamma)$ can be thought of as a function of $\theta\in E$ which describes the exponential growth rate of orbit points converging to limit points of slope $\theta$. It is an invariant of $\Gamma$ which carries finer information than the critical exponent $\delta(\Gamma)$: The critical exponent is simply the maximum of $\delta_\theta(\Gamma)$ among all $\theta\in E$. 
As in  \cite{MR1933790} it will be convenient to study the homogeneous extension  
$$\Psi_\Gamma:\RR_{\ge 0}^r\to \RR$$
of this function $\theta\mapsto \delta_\theta(\Gamma)$. Our final result appropriately
generalizes Theorem~4.2.2 in \cite{MR1933790} (which holds for symmetric spaces and Euclidean buildings of higher rank) to our setting:\\[2mm]
{\bf Theorem E}$\quad$ {\sl $\Psi_\Gamma$ is upper semi-continuous, concave, and the set  
\[ \{L\in \RR_{\ge 0}^r: \Psi_\Gamma(L)>-\infty\}\]
is precisely the limit cone $\ell_\Gamma$ of $\,\Gamma$.
Moreover, $\Psi_\Gamma$ is non-negative on the limit cone $\ell_\Gamma$  and strictly positive on its relative interior. }\\[3mm]
An easy corollary of Theorem~E is the fact that the critical exponent $ \delta(\Gamma)$ of $\Gamma$ is strictly positive and that there exists a unique $\theta^*\in E$ \st $\delta_{\theta^*}(\Gamma)= \delta(\Gamma)$. With the help of 
Theorem~E it is possible 
 to construct generalized conformal densities on each $\Gamma$-invariant subset of the geometric limit set as performed in \cite{MR2062761} and \cite{MR1935549} for higher rank symmetric spaces and Euclidean buildings; this will be done in a future work.
 
The paper is organized as follows: Section~2  summarizes 
basic facts about Hada\-mard spaces  and rank one isometries, in Section~3  specific properties of products of Hadamard spaces are collected. Section~4 provides the main tools needed in the proofs of our results. In Section~5 we describe the structure of the limit set and prove Theorems~A and~B; 
Section~6 is devoted to the study of the limit cone and the proof of Theorem~C. Finally, in Section~7  we introduce and study the homogeneous function $\Psi_\Gamma$ as above 
and prove Theorems~D and~E. \\[-2mm]


\section{Preliminaries}\label{Prelim}

The purpose of this section is to introduce terminology and notation and to summarize basic results about Hadamard spaces and rank one isometries. 
The main references here  are \cite{MR1744486} and  \cite{MR1377265} (see also \cite{MR1383216}, and \cite{MR823981},\cite{MR656659} in the case of Hadamard manifolds). 

Let $(\XX,d)$ be a metric space. A {\hl geodesic path} joining $x\in\XX$ to $y\in\XX$  is a map $\sigma$ from a closed interval $[0,l]\subset \RR$ to $\XX$ \st $\sigma(0)=x$, $\sigma(l)=y$ and $d(\sigma(t), \sigma(t'))=|t-t'|$ for all $t,t'\in [0,l]$.  We will denote such a geodesic path $\sigma_{x,y}$. $\XX$ is called {\hl geodesic} if any two points in $\XX$ can be connected by a geodesic path, if this path is unique, we say that $\XX$ is {\hl uniquely geodesic}. In this text $\XX$ will be a Hadamard space, i.e. a complete geodesic metric space in which all triangles satisfy the CAT$(0)$-inequality. This implies in particular that $\XX$ is simply connected and uniquely geodesic.   A {\hl geodesic} or {\hl geodesic line} in $\XX$ is a map $\sigma:\RR\to\XX$ \st $d(\sigma(t), \sigma(t'))=|t-t'|$ for all $t,t'\in\RR$, a {\hl geodesic ray} is a map $\sigma:[0,\infty)\to \XX$ \st $d(\sigma(t), \sigma(t'))=|t-t'|$ for all $t,t'\in [0,\infty)$. Notice that in the non-Riemannian setting completeness of $\XX$ does not imply that every geodesic path or ray can be extended to a geodesic, i.e. $\XX$ need not be geodesically complete.


From here on we will assume that $\XX$ is a locally compact 
Hadamard space.  The geometric boundary $\rand$ of
$\XX$ is the set of equivalence classes of asymptotic geodesic
rays endowed with the cone topology (see e.g. \cite[chapter~II]{MR1377265}). The action of the isometry group $\is(\XX)$ on $\XX$ naturally extends to an action by homeomorphisms on the geometric boundary. Moreover, since $\XX$ is locally compact, this boundary $\rand$ is compact and the space $\XX$ is a dense and open subset of the compact space $\ganz:=\XX\cup\rand$.  For $x\in\XX$ and $
\xi\in\rand$ arbitrary there exists a  geodesic ray emanating from $x$ which belongs to the class of $\xi$. We will denote such a ray $\sigma_{x,\xi}$.

We say that two points $\xi$, $\eta\in\rand$ can be joined by a geodesic if there exists a geodesic $\sigma:\RR\to\XX$ \st $\sigma(-\infty)=\xi$ and $\sigma(\infty)=\eta$. It is well-known that if all triangles in $\XX$ satisfy the CAT$(-1)$-inequality, 
then every pair of distinct points in the geometric boundary can be joined by a geodesic. This is not true in general. For convenience we therefore define the {\hl visibility set at infinity} $\horinf(\xi)$ of a point  $\xi\in\rand$ as  the set of points in the geometric boundary which can be joined to $\xi$ by a geodesic, i.e.
\begin{equation}\label{visinf}
\horinf(\xi):=\{\eta\in\rand :  \exists\ \mbox{geodesic}\ \sigma\  \st 
\sigma(-\infty)=\xi\,,\,\sigma(\infty)=\eta\}\,.
\end{equation}
%

A geodesic $\sigma: \RR\to\XX$  is said to {\hl bound a flat half-plane}  if there exists a closed convex subset $\iota([0,\infty)\times \RR)$ in $\XX$ isometric to $[0, \infty)\times\RR$ \st $\sigma(t)=\iota(0,t)$ for all $t\in \RR$; 
otherwise $\sigma$ will be called a 
{\hl rank one geodesic}. 

Notice that the existence of a rank one geodesic  imposes severe restrictions on the CAT$(0)$-space $\XX$. For example,  $\XX$ can neither be a symmetric space or  Euclidean building of higher rank  nor a product of Hadamard spaces.

The following important lemma states that even though we cannot join any two distinct points in the geometric boundary of $\XX$, given a rank one geodesic we can at least join points in a neighborhood of its extremities. More precisely, we have the following well-known
\begin{lem}\label{joinrankone} (\cite{MR1377265}, Lemma III.3.1)\ 
Let $\sigma:\RR\to\XX$ be a rank one geodesic. Then there exist $c>0$ and 
neighborhoods $U^-$, $U^+ $ of $\sigma(-\infty)$, 
$\sigma(\infty)$ in $\ganz$ \st for any $\xi\in U^-$ and $\eta \in U^+$ there exists a rank one geodesic joining $\xi$ and $\eta$. For any such geodesic $\sigma'$ we have $d(\sigma', \sigma(0))\le c$.
\end{lem}
Moreover, we will need the following technical lemma which immediately follows from Lemma~4.3 and Lemma~4.4 in \cite{MR1383216}.
\begin{lem}\label{convergence}
Let $\sigma:\RR\to\XX$ be a rank one geodesic and set $y:=\sigma(0)$, $\eta:=\sigma(\infty)$. Then for any $T\gg1$, $\eps>0$  there exists a neighborhood $U^-$ of $\sigma(-\infty)$ in $\ganz$ and a number $R>0$ \st for any $x\in \XX$ with $d(x,\sigma)>R$ or $x\in U^-$ we have
$$ d(\sigma_{x,y}(t),\sigma_{x,\eta}(t))\le \eps\ \quad\mbox{for all}\ \ t\in [0,T].$$
\end{lem}

The following kind of isometries will play a central role in the sequel. 
\begin{df} \label{axialisos}
An isometry $h$ of $\XX$ is called {\hd axial}, if there exists a constant\break
$l=l(h)>0$ and a geodesic $\sigma$ \st
$h(\sigma(t))=\sigma(t+l)$ for all $t\in\RR$. We call
$l(h)$ the {\hd translation length} of $h$, and $\sigma$
an {\hd axis} of $h$. The boundary point
$h^+:=\sigma(\infty)$ is called the {\hd attractive fixed
point}, and $h^-:=\sigma(-\infty)$ the {\hd repulsive fixed
point} of $h$. We further set
$\Ax(h):=\{ x\in\XX:  d(x,h x)=l(h)\}$.
\end{df}
We remark that $\Ax(h)$ consists of the union of parallel geodesics
translated by $h$, and 
$\overline{\Ax(h)}\cap\rand$ is exactly the set of fixed points of
$h$. 
Following the definition in \cite{MR2507218} and \cite{MR2585575} we will call two axial isometries $g$, $h\in\is(\XX)$ {\hl independent} if for any given $x\in\XX$ the map 
$$ \ZZ\times\ZZ\to [0,\infty),\quad (m,n)\mapsto d(g^m   x, h^m  x)$$
is proper.
\begin{df}
An axial isometry  is called {\hd rank one} if it possesses a rank one axis. 
\end{df}
Notice that if $h$ is rank one, then $h^+$ and $h^-$ are the only fixed points of $h$. Moreover, it is easy to verify that two rank one elements $g$, $h\in\is(\XX)$ are independent if and only if $\{g^+,g^-\}\cap\{h^+,h^-\}=\emptyset$. 
Let us recall some properties of rank one isometries.
\begin{lem}\label{dynrankone}(\cite{MR1377265}, Lemma III.3.3)\ 
Let $h$ be a rank one isometry. Then
\begin{enumerate}
\item[(a)] $\horinf(h^+)=\rand\setminus\{h^+\}$, 
\item[(b)] any geodesic joining a point $\xi\in\rand\setminus\{h^+\}$ to $h^+$ is rank one, 
\item[(c)] given neighborhoods $U^-$ of $h^-$ and $U^+$ of $h^+$ in $\ganz$ 
there exists $N_0\in\NN$ \st\\
 $h^{-n}(\ganz\setminus U^+)\subset U^-$ and
$h^{n}(\ganz\setminus U^-)\subset U^+$ for all $n> N_0$.
\end{enumerate}
\end{lem}
The following lemma allows to find many rank one isometries; it will play an important role in the proof of several results such as Proposition
~\ref{constraxial}, which in turn is needed for Proposition~\ref{constrfreesubgrps}.
\begin{lem}\label{elementsarerankone} (\cite{MR1377265}, Lemma III.3.2)\
Let $\sigma:\RR\to\XX$ be a rank one geodesic, and $(\gamma_n)\subset\is(\XX)$ a sequence of isometries \st $\gamma_n x\to \sigma(\infty)$ and $\gamma_n^{-1}x\to\sigma(-\infty)$ for one (and hence any) $x\in\XX$. Then for $n$ sufficiently large, $\gamma_n$ is axial and possesses an axis $\sigma_n$ \st $\sigma_n(\infty)\to\sigma(\infty)$ and $\sigma_n(-\infty)\to\sigma(-\infty)$.
\end{lem}
The following proposition is the clue to the proof of all results in Section~\ref{limsetdiscreteisom} and~\ref{sectlimcone}.
It is more general, but similar in spirit to Proposition~2.8 
in \cite{MR2629900} which was inspired by the proof of Lemma~4.1 in \cite{MR1703039} in the easier context  of CAT$(-1)$-spaces. Part (b) in particular gives a relation between the geometric length and the combinatorial length of words in a free group on two generators which will be the key ingredient in the proof of Proposition~\ref{constrfreesubgrps}. 
If $\alpha,\beta$ generate a free group we say that a word $\gamma=s_1^{k_1}s_2^{k_2}\cdots s_n^{k_n}$ with $s_j\in \{\alpha,\alpha^{-1}, \beta,\beta^{-1}\}$ and $k_j\in\NN$, $j\in\{1,2,\ldots n\}$   is {\hl cyclically reduced} if $s_{j+1}\notin \{s_j, s_j^{-1}\}$ for $j\in\{1,2,\ldots  n-1\}$ and $s_n\ne s_1^{-1}$.  

\begin{prp}\label{combgeomlength}
Suppose  $g$ and $h$ are independent rank one elements in $\is(\XX)$. 
Then there exist neighborhoods $V(\eta), U(\eta)$ of $\eta\in \{g^-,g^{+},h^-,h^{+}\}$ with $V(\eta)\subset U(\eta)\subset\ganz$ and a constant $c>0$ 
\st the following holds:
\begin{enumerate}
\item[(a)]    Any two points in different sets $U(\eta)$ can be joined by a rank one geodesic $\sigma\subset\XX$ with $d(\xo,\sigma)\le c$. 
\item[(b)] Given a pair of rank one isometries $\alpha,\beta\in\is(\XX)$ with $\alpha^\pm\in V(g^\pm)$, $\beta^\pm\in V(h^\pm)$ 
and $N_\alpha, N_\beta\in\NN$ sufficiently large, then for all $n\in\NN$ and every cyclically reduced word
$\gamma=s_1^{k_1}s_2^{k_2}\cdots s_n^{k_n}$ with $s_j\in S:=\{\alpha^{N_\alpha}, \alpha^{-N_\alpha}, \beta^{N_\beta}, \beta^{-N_\beta}\}$ and $k_j\in\NN$,\break  $j\in\{1,2,\ldots n\}$, we have
\begin{enumerate}
\item[(i)]
$\quad \displaystyle\gamma^+\in \Bigg\{\begin{array}{ccc} 
U(g^\pm)&\text{if} & s_1=\alpha^{\pm N_\alpha}\\
U(h^\pm)&\text{if} & s_1=\beta^{\pm N_\beta}\end{array}\quad\text{and}\quad 
\gamma^-\in \Big\{\begin{array}{ccc} 
U(g^\pm)&\text{if} & s_n^{-1}=\alpha^{\pm N_\alpha}\\
U(h^\pm)&\text{if} & s_n^{-1}=\beta^{\pm N_\beta}\end{array}$
\item[(ii)]
$\quad\displaystyle \big|l(\gamma)-\sum_{j=1}^n k_j l(s_j)\big|\le 4c\cdot n\quad\text{and}\quad  \big|d(\xo,\gamma\xo)-\sum_{j=1}^n k_j l(s_j)\big|\le 4c\cdot n$
\end{enumerate}
\end{enumerate}
\end{prp}
\prf\  We fix a base point $\xo\in\Ax(h)$. For $\eta\in \{g^-,g^{+},h^-,h^{+}\}$ let $U(\eta)\subset \ganz$ be a small neighborhood of $\eta$ with $\xo\notin U(\eta)$ \st all $U(\eta)$ are pairwise disjoint, and $c>0$ a constant \st any pair of points  in distinct neighborhoods can be joined by a rank one geodesic $\sigma$ with $d(\xo,\sigma)\le c$. This is possible by Lemma~\ref{joinrankone} and proves part~(a). 

According to Lemma~\ref{elementsarerankone}, for $\eta\in \{g^-,g^{+},h^-,h^{+}\}$ there exist neighborhoods $W(\eta)\subset U(\eta)$ \st every $\gamma\in\Gamma$ with $\gamma\xo\in W(\eta)$, $\gamma^{-1}\xo\in W(\zeta)$, $\zeta\neq \eta$, is rank one with $\gamma^+\in U(\eta)$ and $\gamma^-\in U(\zeta)$.  We claim that assertion (b) holds for all neighborhoods
$V(\eta)\subset W(\eta)\subset\ganz$ of $\eta\in \{g^-,g^{+},h^-,h^{+}\}$ with $c>0$ as above. 

Let $\alpha,\beta\in\is(\XX)$ be rank one isometries with $\alpha^\pm\in W(g^\pm)$,  $\beta^\pm\in W(h^\pm)$ and set $W(\alpha^\pm)= W(g^\pm)$ and   $W(\beta^\pm)= W(h^\pm)$. 
 By Lemma~\ref{dynrankone} (c) there exist $N_\alpha,N_\beta\in\NN$ \st 
\begin{equation}\label{pingpong} 
\alpha^{\pm N_\alpha} \big(\ganz\setminus W(\alpha^\mp)\big)\subset W(\alpha^\pm)\quad\text{and}\quad \beta^{\pm N_\beta} \big(\ganz\setminus W(\beta^\mp)\big)\subset W(\beta^\pm).
\end{equation}
We set $U(\alpha^\pm)=U(g^\pm)$, $U(\beta^\pm)=U(h^\pm)$,  
 $S:=\{\alpha^{N_\alpha},\alpha^{-N_\alpha},\beta^{N_\beta},\beta^{-N_\beta}\}$ and   consider a cyclically reduced word 
$\gamma=s_1^{k_1}s_2^{k_2}\cdots s_n^{k_n}$ with $s_j\in S$ and $k_j\in\NN$, $j\in\{1,2,\ldots n\}$. 
By choice of $N_\alpha,N_\beta$ and (\ref{pingpong}) we have $\gamma\xo\in W(s_1^+)$ and $\gamma^{-1}\xo\in W(s_n^{-})\neq W(s_1^+)$ since $s_1\neq s_n^{-1}$. Therefore $\gamma$ is rank one with $\gamma^+\in U(s_1^+)$ and $\gamma^-\in U(s_n^{-})$, which shows (i). 

Choosing a point $x\in\Ax(\gamma)$ with $d(\xo,x)\le c$ (which is possible according to (a)) we get
\begin{equation}\label{translengthdistance}
 l(\gamma)\le d(\xo,\gamma\xo)\le d(\xo,x)+d(x,\gamma x)+ d(\gamma x,\gamma\xo)\le l(\gamma)+2c\,.
\end{equation}
Similarly, for all $j\in\{1,2,\ldots n\}$ we get $l(s_j^{k_j})\le d(\xo,s_j^{k_j}\xo)\le l(s_j^{k_j})+2c$.

We fix $j\in\{1,2,\ldots n\}$ and abbreviate $\gamma_j:=s_j^{k_j}s_{j+1}^{k_{j+1}}\cdots s_n^{k_n}$. Then $\gamma_2\xo\in W(s_2^+)$, $s_1^{-k_1}\xo\in W(s_1^{-})\neq W(s_2^+)$, so there exists a geodesic $\sigma_2$ joining $\gamma_2\xo$ to $s_1^{-k_1}\xo$ with $d(\xo,\sigma_2)\le c$. If $y$ denotes a point on $\sigma_2$ with $d(\xo,y)\le c$ we get
$$ d(s_1^{k_1}s_2^{k_2}\cdots s_n^{k_n}\xo,\xo)= d(\gamma_2\xo, s_1^{-k_1}\xo)= d(\gamma_2\xo, y)+d(y,s_1^{-k_1}\xo)$$ 
which proves $|d(\gamma\xo,\xo)-d(\xo,s_1^{k_1}\xo)-d(\xo,\gamma_2\xo)|\le 2c$. Applying the same arguments to $\gamma_j$ for $j\ge 2$ and using the fact that $s_{j+1}\ne s_j^{-1}$ we deduce 
\[|d(\gamma_j\xo,\xo)-d(\xo,s_j^{k_j}\xo)-d(\xo,\gamma_{j+1}\xo)|\le 2c\] and hence inductively
\[ \big| d(\xo,\gamma\xo)-\sum_{j=1}^n d(\xo,s_j^{k_j}\xo)\big|\le 2(n-1)c. \]
Using (\ref{translengthdistance})  we conclude
\[ \big|l(\gamma)-\sum_{j=1}^n k_j l(s_j)\big|\le 4c\cdot n\quad\text{and}\quad 
\big|d(\xo,\gamma\xo)-\sum_{j=1}^n k_j l(s_j)\big|\le 4c\cdot n.\]
Therefore  part~(ii)  of (b) is also true.\qed 

\section{Products of Hadamard spaces}\label{prodHadspaces}

Now let $(\XX_1,d_1)$, $(\XX_2,d_2),\ldots, (\XX_r,d_r)$ be locally compact 
Hadamard spaces,  and $\XX=\XX_1\times \XX_2\times \cdots\times \XX_r$  the product space endowed with the product distance $d=\sqrt{d_1^2+d_2^2+\cdots + d_r^2}$. Notice that such a product is again a locally compact 
Hadamard space.  

We denote $\RR_{\ge 0}^r:=\big\{(t_1,t_2,\ldots,t_r)\in\RR^r: t_i\ge 0\ \text{for all}\ i\in\{1,2,\ldots,r\}\big\}$ and\break
$\RR_{>0}^r:=\big\{(t_1,t_2,\ldots,t_r)\in\RR^r: t_i > 0\ \text{for all}\ i\in\{1,2,\ldots,r\}\big\}$. 
To any pair of points $x=(x_1,x_2,\ldots, x_r)$, $z=(z_1,z_2,\ldots, z_r)\in \XX$ we associate the vector 
$$H(x,z):= \left(\begin{array}{c} d_1(x_1,z_1)\\ d_2(x_2,z_2)\\\vdots\\ d_r(x_r,z_r)\end{array}\right)\in \RR_{\ge 0}^r\,,$$ 
which we call the {\hl distance vector} of the pair $(x,z)$. Notice that  if $\Vert \cdot\Vert$ denotes the Euclidean norm in $\RR^r$, we clearly have $\Vert H(x,z)\Vert =d(x,z)$. If $z\neq x$ we therefore define the {\hl direction} of $z$  \wrt $x$ by the unit vector 
$$\thet(x,z):=\frac{H(x,z)}{d(x,z)}  \in\RR_{\ge 0}^{r}.$$

Denote $p_i:\XX\to \XX_i$, $i=1,2,\ldots, r$, the natural projections. Every 
geodesic path $\sigma:[0,l]\to\XX\,$ can be written as a product $\sigma(t)=(\sigma_1(t\cdot \theta_1), \sigma_2(t\cdot \theta_2),\ldots, \sigma_r(t \cdot \theta_r))$, where  $\sigma_i$ are geodesic paths in $\XX_i$, $i=1,2,\ldots r$, and the $\theta_i\ge 0$ satisfy 
\[\sum_{i=1}^r \theta_i^2=1.\]  
The unit vector
\[ \text{sl}(\sigma):= \left(\begin{array}{c} \theta_1\\ \theta_2\\\vdots\\ \theta_r\end{array}\right)\in E:=\{\theta\in  \RR_{\ge 0}^r :\ \Vert \theta\Vert =1\}\]
equals the direction of the points $\sigma(t)$, $t\in (0,l]$, with respect to $\sigma(0)$ and  is called the {\hl slope of $\sigma$}.  We say that a geodesic path $\sigma$ is  {\hl regular} if its slope does not possess a coordinate zero, i.e. if 
\[ \text{sl}(\sigma)\in E^+:=\{\theta\in  \RR_{> 0}^r :\ \Vert \theta\Vert =1\};\]
otherwise $\sigma$ is said to be {\hl singular}.  In other words, $\sigma$ is regular if and only if none of the projections $p_i(\sigma([0,l]))$, $i\in\{1,2,\ldots, r\}$,   is a point.

If  $x\in\XX$  and $\sigma:[0,\infty)\to\XX$ is an arbitrary geodesic ray, then elementary geometric estimates yield  the relation 
\[ 
 \text{sl}(\sigma)= \lim_{t\to\infty} \thet(x,\sigma(t)) =\lim_{t\to\infty} \frac{H(x,\sigma(t))}{d(x,\sigma(t))}
\]
 between the slope of $\sigma$ and the directions of $\sigma(t)$, $t>0$, \wrt $x$. 
 Similarly, one can easily show that any two geodesic rays representing the same (possibly singular) point in the geometric boundary necessarily have the same slope. So  the slope $\text{sl}(\tilde\xi)\in E$ of a point $\tilde\xi\in\rand$ can be defined as the slope of an arbitrary geodesic ray representing $\tilde\xi$.  The {\hl regular geometric boundary} $\regrand$ and  the {\hl singular geometric boundary} $\singrand$ of $\XX$ are then naturally defined by 
 \[ \regrand:=\{\tilde\xi\in\rand:\ \text{sl}(\tilde\xi)\in E^+\},\quad \singrand:=\rand\setminus\singrand;\]
the singular boundary $\singrand$ consists of equivalence classes of geodesic rays in $\XX$ which project to a point in one of the factors $\XX_i$.

We further notice that two regular geodesic rays  $\sigma$, $\sigma'$  with the same slope represent the same point in $\regrand$ if and only if  $\sigma_i(\infty)=\sigma_i'(\infty)$ for all $i\in\{1,2,\ldots,r\}$.  So $\regrand$ is homeomorphic
to $\rand_1\times\rand_2\times\cdots\times\rand_r\times E^+$.

 If $\gamma\in\is(\XX_1)\times \is(\XX_2)\times\cdots\times\is(\XX_r)$, then 
 the slope of $\gamma\at\tilde\xi$ equals the slope of $\tilde\xi$. In other words, if $\rand_\theta$ denotes the set of points in the geometric boundary of slope $\theta\in E $, then $\rand_\theta$ is invariant by the action of $\is(\XX_1)\times\is(\XX_2)\times\cdots\times\is(\XX_r)$.  

In analogy to the case of symmetric spaces of higher rank we define the {\hl Furstenberg boundary} $\Frand$ of $\XX$ as the product $\rand_1\times \rand_2\times\cdots\times \rand_r$ endowed with the product topology. Since $\regrand$ is homeomorphic to $\Frand\times E^+$ 
we have a natural projection
$$ \hspace{1cm} \ba{lccc}\pi^F\,:& \regrand&\to& \Frand\\
&(\xi_1,\xi_2,\ldots,\xi_r,\theta)&\mapsto & (\xi_1,\xi_2,\ldots,\xi_r)\, \ea $$
and a natural action of the  group $\is(\XX_1)\times\is(\XX_2)\times\cdots\times \is(\XX_r)$ by homeomorphisms on the Furstenberg boundary of $\XX$. 
Clearly $\Frand$ is homeomorphic to each of the sets $\rand_\theta\subset\regrand$ with $\theta\in E^+$. Notice that in the special case $r=1$ the Furstenberg boundary $\Frand$ equals the geometric boundary and the projection $\pi^F$ is the identity; $E^+=E$ is simply a point. 


The following two elementary lemmata provide important facts concerning the topology of $\ganz$. 
\begin{lem}\label{prodtopologytheta}
Suppose $(y_n)\subset\XX$ is a sequence converging to a point $\tilde\eta\in\rand_{\theta}$ for some $\theta \in E$. Then for any $x\in\XX\,$ we have $\thet(x,y_n)\to \theta$ as $n\to\infty$. 
\end{lem}
Notice that if $\theta=(\theta_1,\theta_2,\ldots, \theta_r)\in E$ satisfies $\theta_i=0$ for some $i\in\{1,2,\ldots,r\}$, then the projections $y_{n,i}=p_i(y_n)$ of $y_n$ to $\XX_i$ necessarily satisfy
\[ \lim_{n\to\infty} \frac{d_i(x_i,y_{n,i})}{d(x,y_n)} =0.\]  
So the sequence $(y_{n,i})\subset \XX_i$ can be either bounded or unbounded and may possess more than one accumulation point in $\ganz_i$.  However, if $\theta\in E^+$ then we have the following
\begin{lem}\label{prodtopologyreg}
Suppose $(y_n)\subset\XX$ is a sequence converging to a regular boundary point $\tilde\eta\in\regrand$ with $\pi^F(\tilde\eta)=(\eta_1,\eta_2,\ldots,\eta_r)\in\Frand$. Then for all $i\in\{1,2,\ldots,r\}$ the projections $p_i(y_n)$ to $\XX_i$ converge to $\eta_i$ as $n\to\infty$. 
\end{lem}

Recall the definition of the visibility set at infinity  $\horinf(\tilde\xi)$  of a point $\tilde\xi\in\rand$ from~(\ref{visinf}). It is easy to see that a point $\tilde\eta\in\rand$ cannot belong to $\horinf(\tilde\xi)$ if the slope of $\tilde\eta$ is different from the slope of $\tilde\xi$. This motivates the less restrictive definition of  the {\hl Furstenberg visibility set} of a point $\xi\in\Frand$ which is
\[ 
 \horF(\xi):=\pi^F\big(\horinf(\tilde\xi)\big)\,,\quad \mbox{where}\ \  \tilde\xi\in (\pi^F)^{-1}(\xi)\ \mbox{is arbitrary}.\]
We say that $\xi,\eta\in\Frand$ are {\hl opposite} if $\eta\in\horF(\xi)$. Notice that $\xi=(\xi_1,\xi_2,\ldots, \xi_r)$, $\eta=(\eta_1,\eta_2,\ldots, \eta_r)\in\Frand$ are opposite if and only if $\xi_i$ and $\eta_i$ can be joined by a geodesic in $\XX_i$ for all $i\in\{1,2,\ldots, r\}$, i.e.
\begin{equation}\label{visFequiv}
\horF(\xi)=\big\{(\eta_1,\eta_2,\ldots,\eta_r)\in\Frand:\eta_i\in\horinf(\xi_i)\ \text{for all}\ i\in\{1,2,\ldots, r\}\big\} . 
\end{equation}

We terminate the section with a few definitions concerning isometries of products.  By abuse of notation we also denote $p_i:\is(\XX_1)\times\is(\XX_2)\times\cdots\times\is(\XX_r)\to \is(\XX_i)$, $i=1,2,\ldots,r$, the natural projections.
\begin{df}\label{regaxial}
An isometry $h=(h_1,h_2,\ldots, h_r)\in\is(\XX_1)\times\is(\XX_2)\times\cdots\times\is(\XX_r)$ is called {\hd regular axial}, if  $h_i=p_i(h)$ is a rank one isometry of $\XX_i$ for all $i\in\{1,2,\ldots,r\}$; we denote $\widetilde{h^+}$ its attractive fixed point in $\regrand$ and set $h^+:=\pi^F(\widetilde{h^+})=(h_1^+, h_2^+,\ldots,h_r^+)$. Moreover,
we denote by $l_i(h)$, $i=1,2,\ldots, r$, the {\hd translation length} of the rank one isometry $p_i(h)$ in $\XX_i$, and by 
\[ L(h):=\left(\begin{array}{c} l_1(h)\\ l_2(h)\\\vdots\\ l_r(h)\end{array}\right)\in\RR_{>0}^r\]
the {\hd translation vector} of $h$.
\end{df}
Notice that (\ref{visFequiv}) and Lemma~\ref{dynrankone} (a)  imply  
\[\horF(h^+)=\big\{(\xi_1,\xi_2,\ldots,\xi_r)\in\Frand : \xi_i\ne h_i^+,\ \text{for all}\ i\in\{1,2,\ldots,r\}\big\}.\]
Moreover, the translation length $l(h)$ in $\XX$ of a regular axial isometry is given by
$l(h)=\Vert L(h)\Vert$; the unit vector 
\[ \widehat L(h):=\frac{L(h)}{l(h)}\in E^+\]
is sometimes called the {\hl translation direction} of $h$.
\begin{df}
Two regular axial isometries $h, g\in\is(\XX_1)\times\is(\XX_2)\times\cdots\times\is(\XX_r)$ are called {\hd independent}, if $p_i(h)$ and $p_i(g)$ are independent for all $i\in\{1,2,\ldots,r\}$.
\end{df}
In other words, $h$ and $g$ are independent if 
\[ \{g^+,g^-\}\subset \horF(h^+)\cap \horF(h^-);\]
such pairs of independent regular axial isometries will play a key role throughout the article.
%

\section{Key results on pairs of independent regular axials}

Recall that $\XX=\XX_1\times \XX_2\times \cdots\times \XX_r$  is a product of locally compact Hadamard spaces. 
We fix a regular axial isometry 
$h=(h_1,h_2,\ldots, h_r)\in\is(\XX_1)\times\is(\XX_2)\times\cdots\times\is(\XX_r)$ and a base point $\xo=(\xo_1,\xo_2,\ldots,\xo_r)\in\Ax(h)$; in particular, for each $i\in\{1,2,\ldots, r\}$ the point $\xo_i\in\XX_i$ lies on an invariant geodesic of the rank one isometry $h_i\in\is(\XX_i)$.  

The results in this section are key ingredients for the proofs of all results. 
Proposition~\ref{movingpointsnbhds}, Lemma~\ref{separatepoints}, Proposition~\ref{gettingnicesequences} and Corollary~\ref{gettingregularaxialsinGamma} generalize Proposition~5.1 and its consequences in~\cite{MR2629900}, where the analogous results in the case of two factors were proved  by a case-by-case study. When considering more factors this method does not work any more, because the projections of points in $\ganz_1\times \ganz_2\times\cdots\times\ganz_r$ to all factors need to be controlled simultanously. 
The new idea for Proposition~\ref{movingpointsnbhds} here is to move the first projection to the desired set without worrying about what happens to the other projections, and then step by step take care of the remaining projections without messing up 
what was already moved to the right place. This idea 
can also be used to show that {\hl two} points in the Furstenberg boundary can be moved by the {\hl same} isometry to an arbitrarily small neighborhood of $h^+$; this is the content of Proposition~\ref{contract2points}. 

In Proposition~\ref{constraxial} we further provide the analogon of Proposition~2.2.7 in~\cite{MR1933790} (compare also \cite{MR1437472}), which turns out to be indispensable for the construction of free semi-groups in a discrete group as performed in Proposition~\ref{constrfreesubgrps}. This construction in turn plays a central role in the proofs of Proposition~\ref{limitconeconvex}, Proposition~\ref{constrfreegrps} and finally Theorem~\ref{expgrowthpositive}. Finally, Proposition~\ref{genprod} states the equivalent of Proposition~2.3.1 in \cite{MR1933790} in our setting, which is necessary for the proof of Theorem~\ref{concave}. Thanks to Proposition~\ref{movingpointsnbhds}, the proof of the corresponding  
Proposition~7.2 in \cite{MR2629900} 
easily extends  to more than two factors.
\begin{prp}\label{movingpointsnbhds}
Assume that $g\in\is(\XX_1)\times\is(\XX_2)\times\cdots\times\is(\XX_r)$ is regular axial with $\{g^+,g^-\}\subset\horF(h^-)$ and $g^-\in\horF(h^+)$. Given neighborhoods $U_i\subset\ganz_i$ of $h_i^+$, $i=1,2,\ldots,r$, there exist $N\in\NN$ and a finite set
$\Lambda_r\subset\langle g^N,h^N\rangle^+$ in the semi-group generated by $g^N$ and $h^N$ consisting of $2^r$ words of length at most $2r$ in the generators $g^N$, $h^N$ 
\st for any $z=(z_1,z_2,\ldots,z_r)\in \ganz_1\times\ganz_2\times\cdots\times \ganz_r$
there exists $\lambda(z)=(\lambda_1,\lambda_2,\ldots,\lambda_r)\in\Lambda_r$ with 
\[\lambda(z)\cdot z\in U_1\times U_2\times\cdots\times U_r.\]
Moreover, if $z=(z_1,z_2,\ldots,z_r)\in \ganz_1\times\ganz_2\times\cdots\times \ganz_r$ 
and $\lambda(z)=(\lambda_1,\lambda_2,\ldots,\lambda_r)\in\Lambda_r$, then there exist open neighborboods $V_i\subset\ganz_i$ of $z_i$, $i=1,2,\ldots r$,    with the property
\[ \lambda_i\cdot V_i\subset U_i\quad\text{for all}\quad i\in \{1,2,\ldots,r\}.\]
\end{prp}
\prf\ For $i=1,2,\ldots, r$ and $\eta_i\in\rand_i$ a point in the set $\{g_i^-,g_i^{+},h_i^-,h_i^{+}\}$ (which consists of only 3 points if $g_i^+=h_i^+$)  let $W_i(\eta_i)\subset\ganz_i$ be an arbitrary sufficiently small neighborhood of $\eta_i^+\in\rand_i$ with $\xo_i\notin W_i(\eta_i)$ \st the closures of 
all $W_i(\eta_i)$ are pairwise disjoint in $\ganz_i$. Making $W_i(h_i^+)$ smaller if necessary we may further assume that 
$W_i(h_i^+)\subset U_i$ for all $i\in\{1,2,\ldots,r\}$. 
 According to Lemma~\ref{dynrankone} (c) there exists  a constant   $N\in\NN$ \st for all  $i\in\{1,2,\ldots,r\}$
\begin{equation}\label{factorpingpong}  
g_i^{\pm N}\big(\ganz_i\setminus  W_i(g_i^{\mp})\big)\subset W_i(g_i^\pm),\qquad h_i^{\pm N}\big(\ganz_i\setminus  W_i(h_i^{\mp})\big)\subset W_i(h_i^\pm).
\end{equation}
We prove the first claim by induction on $r$: 

For $r=1$ we let $z=z_1\in\ganz_1=W_1(h_1^-)\cup \ganz_1\setminus \overline{W_1(h_1^-)}$.  If $z_1\in  W_1(h_1^-)$, then from  $W_1(h_1^-)\subset \ganz_1\setminus W_1(g_1^-)$ and~(\ref{factorpingpong}) we get $g_1^N z_1\in W_1(g_1^+)\subset \ganz_1\setminus  W_1(h_1^-)$, hence again by ~(\ref{factorpingpong})
\[h_1^N g_1^N z_1\in W_1(h_1^+)\subset U_1.\]
If $z_1\in\ganz_1\setminus \overline{ W_1(h_1^-)}\subset \ganz_1\setminus W_1(h_1^-) $, then~(\ref{factorpingpong}) directly gives $h_1^Nz_1\in W_1(h_1^+)\subset U_1$. 
In particular this shows that $\Lambda_1:=\{h^Ng^N, h^N\}\subset\langle g^N,h^N\rangle^+$ is the desired set consisting of $2=2^1$ elements of length $\le 2=2\cdot 1$.

Moreover, since both $W_1(h_1^-)$ and $\ganz_1\setminus \overline{W_1(h_1^-)}$ are open, there exists an open neighborhood $V_1\subset\ganz_1$ of $z_1$ \st either
\[ h_1^N g_1^N \cdot V_1\subset U_1\qquad\text{or}\quad h_1^N\cdot V_1\subset U_1.\]


Now assume the assertion holds for $r-1$;  we claim that it also holds when $r$ factors are involved.  By the induction hypothesis there exists a finite set \[\Lambda_{r-1}\subset \langle (g_1,g_2,\ldots, g_{r-1})^N,(h_1,h_2,\ldots, h_{r-1})^N\rangle^+<\is(\XX_1)\times\is(\XX_2)\times\cdots\times\is(\XX_{r-1})\] 
consisting of $2^{r-1}$ words of length at most $2(r-1)$ in its generators \st 
for all $(y_1,y_2,\ldots, y_{r-1})\in  \ganz_1\times \ganz_2\times\cdots\times \ganz_{r-1}$
there exists $\lambda'=(\lambda_1',\lambda_2',\ldots,\lambda_{r-1}')\in\Lambda_{r-1}$ 
and open neighborhoods $V_i\subset\ganz_i$ of $y_i$   
with $\lambda_i'  \cdot V_i \subset U_i$ for all $i\in\{1,2,\ldots,r-1\}$.
We denote by  $\Lambda'\subset \langle g^N,h^N\rangle^+$ the finite set of the same words as in $\Lambda_{r-1}$, but now considered as elements in $\is(\XX_1)\times\is(\XX_2)\times\cdots\times\is(\XX_{r})$, and let $z=(z_1,z_2,\ldots, z_r)\in \ganz_1\times \ganz_2\times\cdots\times \ganz_r$ arbitrary. 
By the properties of $\Lambda_{r-1}$ we know that there exists $\lambda'=(\lambda_1',\lambda_2',\ldots,\lambda_{r-1}',\lambda_r') \in\Lambda'$ and open neighborhoods $V_i\subset\ganz_i$ of $z_i$, $1\le i\le r-1$    
\st  
\[ \lambda_i'  \cdot V_i \subset U_i\quad\text{for all} \quad i\in\{1,2,\ldots,  r-1\},\]
but we do not know the position of $\lambda_r' z_r\in\ganz_r=W_r(h_r^-)\cup \ganz_r\setminus \overline{W_r(h_r^-)}$. 

However, as in the case $r=1$ the north-south dynamics~(\ref{factorpingpong}) implies 
\[ h_r^N g_r^N \lambda_r'z_r\in U_r\qquad\text{or}\quad h_r^N \lambda_r'z_r\in U_r\]
according to the two cases $\lambda_r'z_r\in W_r(h_r^-)$ or $\lambda_r' z_r\in\ganz_r\setminus \overline{W_r(h_r^-)}$. Moreover, there exists an open neighborhood $V_r\subset\ganz_r$ of $z_r$ \st $\lambda_r'\cdot V_r\subset W_r(h_r^-)$ or $\lambda_r' \cdot V_r\subset\ganz_r\setminus \overline{W_r(h_r^-)}$ and hence
\[ h_r^N g_r^N \lambda_r'\cdot V_r\subset U_r\qquad\text{or}\quad h_r^N \lambda_r'\cdot V_r\subset U_r.\]
%
Since for all $i\in\{1,2,\ldots, r-1\}$ we have
\[ h_i^N g_i^N\cdot  U_i\subset W_i(h_i^+)\subset U_i\quad\text{and}\quad h_i^N \cdot U_i\subset  U_i\] we conclude that 
 the set $\Lambda_r$ consisting of all words in $g^N, h^N$ of the form
$h^N\lambda'$ or $h^Ng^N\lambda'$ with $\lambda'\in\Lambda'$ works. Clearly, all such words have length $\le 2 + 2(r-1)=2r$ in the generators $g^N, h^N$ and the cardinality of $\Lambda_r$ is equal to $2\cdot 2^{r-1}=2^r$.\qed\\[3mm]
%
\begin{rmke} If $g$ and $h$ are independent, then -- replacing $h$ by $h^{-1}$ -- an analogous statement holds for neighborhoods $U_i$ of $h_i^-$, $i=1,2,\ldots, r$. 
\end{rmke}

In this case we have the following easy corollary which will be used in the proof of Theorem~\ref{boundlimset}. 
\begin{lem}\label{separatepoints}
Assume that $g\in\is(\XX_1)\times\is(\XX_2)\times\cdots\times\is(\XX_r)$ is regular axial and $g$, $h$ are independent. Then for any $\zeta\in\Frand$ and all $\eta\in\Flim$ there exists $\alpha\in\Gamma$ \st
$\alpha\zeta\in\horF(\eta)$.
\end{lem}
\prf \  For $i\in\{1,2,\ldots,r\}$ and $\eta_i\in \{h_i^{+},h_i^{-}\}$ we let $U_i(\eta_i)\subset\ganz_i$ be the neighborhoods satisfying property (a) of Proposition~\ref{combgeomlength}. According to Proposition~\ref{movingpointsnbhds} there exist $\lambda,\mu\in\Gamma$ \st 
\[\lambda\zeta\in U_1(h_1^+)\times U_2(h_2^+)\times\cdots\times U_r(h_r^+)\quad\text{and}\quad \mu\eta\in U_1(h_1^-)\times U_2(h_2^-)\times\cdots\times U_r(h_r^-),\]
which immediately gives $\lambda\zeta\in\horF(\mu\eta)$ and hence $\mu^{-1}\lambda\zeta\in\horF(\eta)$.\qed\\[3mm]
We next state a stronger version of Proposition~\ref{movingpointsnbhds} which will also  be needed in the proof of Theorem~\ref{boundlimset}. 
\begin{prp}\label{contract2points}
Assume that $g\in\is(\XX_1)\times\is(\XX_2)\times\cdots\times\is(\XX_r)$ is regular axial and\break $g$, $h$ are independent. Then for any neighborhood $U\subset\Frand$ of $h^+$ there exists a finite set $\Lambda\subset\langle g,h\rangle$ \st for any two points $\zeta,\eta\in\Frand$ and some $\lambda\in\Lambda $ we have $\{\lambda\zeta,\lambda\eta\}\subset U$.
\end{prp}
\prf\ For $i=1,2,\ldots, r$ and $\eta_i\in\rand_i$ a point in the set $\{g_i^-,g_i^{+},h_i^-,h_i^{+}\}$ we  let $W_i(\eta_i)\subset\ganz_i$ be an arbitrary 
neighborhood of $\eta_i^+\in\rand_i$ with $\xo_i\notin W_i(\eta_i)$ \st the closures of 
all $W_i(\eta_i)$ are pairwise disjoint in $\ganz_i$. Making $W_i(h_i^+)$ smaller if necessary we may further assume that 
$W_1(h_1^+)\times W_2(h_2^+)\times\cdots\times W_r(h_r^+)\subset U$.
 According to Lemma~\ref{dynrankone} (c) there exists  a constant   $N\in\NN$ \st for all  $i\in\{1,2,\ldots,r\}$
\[g_i^{\pm N}\big(\ganz_i\setminus  W_i(g_i^{\mp})\big)\subset W_i(g_i^\pm),\qquad h_i^{\pm N}\big(\ganz_i\setminus  W_i(h_i^{\mp})\big)\subset W_i(h_i^\pm).\]
We  prove the claim by induction on $r$:\  
For $r=1$ we let $\zeta=\zeta_1\in\rand_1$ and $\eta=\eta_1\in\rand_1$ arbitrary. If both $\zeta$ and $\eta$ belong to $W_1(h_1^-)$, then $\lambda= h^Ng^N$ is the desired element, if both $\zeta$ and $\eta$ are contained in $\ganz_1\setminus  
W_1(h_1^-)$, then $\lambda= h^N$ is. So it remains to deal with the case that one of the points, say $\zeta$, belongs to $W_1(h_1^-)$ and the second one does not: If $\eta\in\ganz_1\setminus W_1(g_1^-)$, then $\lambda= h^Ng^N$ works again, if $\eta\in W_1(g_1^-)$ we can take $\lambda= h^Ng^{-N}$. In particular there exists $\lambda\in\Lambda:=\{h^N, h^Ng^N, h^Ng^{-N}\}$ \st $\{\lambda\zeta,\lambda\eta\}\subset W_1(h_1^+)\subset U$.

Now assume that there exists a finite set \[\Lambda_{r-1}\subset \langle (g_1,g_2,\ldots, g_{r-1}),(h_1,h_2,\ldots, h_{r-1})\rangle\subset\is(\XX_1)\times\is(\XX_2)\times\cdots\times\is(\XX_{r-1})\] 
 \st for any 
$(\zeta_1,\zeta_2,\ldots, \zeta_{r-1})$, $(\eta_1,\eta_2,\ldots, \eta_{r-1})\in  \rand_1\times \rand_2\times\cdots\times \rand_{r-1}$
and some $\lambda'=(\lambda_1',\lambda_2',\ldots,\lambda_{r-1}')\in\Lambda_{r-1}$ we have 
\[\{\lambda_i'\zeta_i, \lambda_i'\eta_i\} \subset W_i(h_i^+)\qquad\text{for all}\quad i\in\{1,2,\ldots,r-1\}.\]
  We denote by  $\Lambda'\subset \langle g,h\rangle$ the finite set of the same words as in $\Lambda_{r-1}$, but now considered as elements in $\is(\XX_1)\times\is(\XX_2)\times\cdots\times\is(\XX_{r})$, and let $\zeta=(\zeta_1,\zeta_2,\ldots, \zeta_{r-1},\zeta_r)$, \break $\eta=(\eta_1,\eta_2,\ldots, \eta_{r-1},\eta_r)\in  \rand_1\times \rand_2\times\cdots\times \rand_{r-1}\times \rand_r=\Frand$ arbitrary. 
By the properties of $\Lambda_{r-1}$ we know that there exists $\lambda'=(\lambda_1',\lambda_2',\ldots,\lambda_{r-1}',\lambda_r') \in\Lambda'$ \st for $i\in\{1,2,\ldots,r-1\}$ we have
$\{\lambda_i'\zeta_i, \lambda_i'\eta_i\} \subset W_i(h_i^+)$. For 
$\lambda_r'\zeta_r$ and $\lambda_r'\eta_r$ there are different possibilities: If both 
points belong to $W_r(h_r^-)$, then $\lambda:= h^N g^N\lambda'$ moves both $\zeta$ and $\eta$ to $W_1(h_1^+)\times W_2(h_2^+)\times\cdots\times W_r(h_r^+)\subset U$, if both points are contained in $\ganz_r\setminus W_r(h_r^-)$, then $\lambda:= h^N\lambda'$ satisfies $\{\lambda\zeta, \lambda\eta\}\subset W_1(h_1^+)\times W_2(h_2^+)\times\cdots\times W_r(h_r^+)\subset U$. It finally remains to deal with the case that one of the points, say $\lambda_r'\zeta_r$, belongs to $W_r(h_r^-)$ and the second one does not: If $\lambda_r'\eta_r\in\ganz_r\setminus W_r(g_r^-)$, then $\lambda= h^Ng^N\lambda'$ works again, if $\lambda_r'\eta_r\in W_r(g_r^-)$ we can take $\lambda= h^Ng^{-N}\lambda'$. So we conclude that for some 
\[\lambda\in\Lambda:=h^N \Lambda'\cup h^Ng^N\Lambda'\cup h^Ng^{-N}\Lambda'\]
 we have  $\quad\{\lambda\zeta,\lambda\eta\}\in W_1(h_1^+)\times W_2(h_2^+)\times\cdots\times W_r(h_r^+)\subset U$.\qed 
\begin{prp}\label{gettingnicesequences}
Let $\Gamma<\is(\XX_1)\times\is(\XX_2)\times\cdots\times\is(\XX_r)$ be a discrete group containing $h$ and a second regular axial element $g$ \st  
$g$ and $h$ are independent.  Let\break 
$(\gamma_n)=\big((\gamma_{n,1},\gamma_{n,2},\ldots, \gamma_{n,r})\big) \subset\Gamma$ be a sequence \st  for all $i\in\{1,2,\ldots,r\}$ the sequences $\gamma_{n,i}\xo_i$ and $\gamma_{n,i}^{-1}\xo_i$ converge to points in $\rand_i$ as $n\to\infty$. Then given arbitrarily small distinct neighborhoods $W_i(h_i^+), W_i(h_i^-)\subset\ganz_i$ of $h_i^+$, $h_i^-$, $i=1,2,\ldots,r$, there exist $N\in\NN$, finite sets $\Lambda^+\subset \langle h^N,g^N\rangle^+$, $\Lambda^-\subset \langle h^{-N},g^{-N}\rangle^+$, $\lambda\in\Lambda^+$  and $\mu\in\Lambda^{-}$ \st $\varphi_n:=\lambda\gamma_n\mu^{-1}$ satisfies 
\begin{align*}
&\varphi_{n}\xo\in W_1(h_1^+)\times W_2(h_2^+)\times\cdots\times W_r(h_r^+)\quad\text{and}\\
 &\varphi_{n}^{-1}\xo\in W_1(h_1^-)\times W_2(h_2^-)\times\cdots\times W_r(h_r^-)\end{align*} for $n$ sufficiently large.
\end{prp}
\prf  For the neighborhoods  $W_i(h_i^+), W_i(h_i^-)\subset\ganz_i$ of $h_i^+$, $h_i^-$, $i\in\{1,2,\ldots,r\}$, we let $N\in\NN$ and $\Lambda^+\subset\langle g^N, h^N\rangle^+$, 
$\Lambda^-\subset\langle g^{-N}, h^{-N}\rangle^+$ be the finite sets according to Proposition~\ref{movingpointsnbhds}. That is for any $z=(z_1,z_2,\ldots,z_r)\in \ganz_1\times\ganz_2\times\cdots\times \ganz_r$ 
there exists $\lambda(z)\in\Lambda^+$ and $\mu(z)\in\Lambda^-$ \st
\begin{align*}
& \lambda(z)\cdot z\in W_1(h_1^+)\times W_2(h_2^+)\times\cdots\times W_r(h_r^+)\quad\text{and}\\
\quad &\mu(z)\cdot z\in W_1(h_1^-)\times W_2(h_2^-)\times\cdots\times W_r(h_r^-).
\end{align*} 

We denote $F\subset\XX$ the finite set of points $\{\lambda^{-1}\xo :\lambda\in\Lambda^{+}\cup\Lambda^-\}\subset\XX$. 
For $i\in\{1,2,\ldots,r\}$ we let $\xi_i\in\rand_i$ be the limit of the sequence $(\gamma_{n,i}\xo_i)\subset\XX_i$. 
By Proposition~\ref{movingpointsnbhds} there exist $\lambda=(\lambda_1,\lambda_2,\ldots,\lambda_r)\in\Lambda^+$ and neighborhoods $V_i^+$ of $\xi_i$ in $\ganz_i$ with $\lambda_i \cdot V_i^+\subset W_i(h_i^+)$ for all $i\in\{1,2,\ldots,r\}$. Since for any $x_i\in\XX_i$ the sequence $\gamma_{n,i}x_i$ also converges to $\xi_i$, there  
exists $N_+\in\NN$ \st for all $n>N_+$ and every $x\in F$ we have $\gamma_{n} x\in V_1^+\times V_2^+\times\cdots\times V_r^+$, and hence 
\[\lambda\gamma_n x \in W_1(h_1^+)\times W_2(h_2^+)\times\cdots\times W_r(h_r^+).\]
Similarly, if $\zeta_i\in\rand_i$ is the limit of the sequence $(\gamma_{n,i}^{-1}\xo_i)\subset\XX_i$, $i\in\{1,2,\ldots,r\}$, 
then 
there exist $\mu=(\mu_1,\mu_2,\ldots,\mu_r)\in\Lambda^-$ and neighborhoods $V_i^-$ of $\zeta_i$ in $\ganz_i$ with $\mu_i \cdot V_i^-\subset W_i(h_i^-)$ for all $i\in\{1,2,\ldots,r\}$. As before, there exists $N_-\in\NN$ \st for all $n>N_-$ and every $x\in F$ we have $\gamma_{n}^{-1} x\in V_1^-\times V_2^-\times\cdots\times V_r^-$, and hence 
\[\mu\gamma_n^{-1} x \in W_1(h_1^-)\times W_2(h_2^-)\times\cdots\times W_r(h_r^-).\]
Since both $\lambda^{-1}\xo$ and $\mu^{-1}\xo$ belong to  the finite set $F$ the assertion is true for all\break  
$n>\max\{N_+,N_-\}$.\qed\\
\begin{rmk} The assumption concerning the sequence $(\gamma_n)$ in $\Gamma$ is clearly satisfied if $\gamma_n\xo$ and $\gamma_n^{-1}\xo$ converge to points in the regular boundary $\regrand$ of $\XX$. However, the result is also valid if 
$\gamma_n\xo$ and $\gamma_n^{-1}\xo$ converge to 
singular boundary points in a way that for all $i\in\{1,2,\ldots,r\}$ the 
sequences  $\gamma_{n,i}\xo_i$ 
and $\gamma_{n,i}^{-1}\xo_i$ 
converge to points in $\rand_i$. 
\end{rmk}\\[3mm]
In combination 
with Lemma~\ref{elementsarerankone}, we get the following useful
\begin{cor}\label{gettingregularaxialsinGamma}
Let $\,\Gamma<\is(\XX_1)\times\is(\XX_2)\times\cdots\times\is(\XX_r)$ be a discrete group containing $h$ and a second regular axial element $g$ \st  $g$ and $h$ are independent.  Let\break
 $(\gamma_n)=\big((\gamma_{n,1},\gamma_{n,2},\ldots, \gamma_{n,r})\big) \subset\Gamma$ be a sequence \st  for all $i\in\{1,2,\ldots,r\}$ the sequences $\gamma_{n,i}\xo_i$ and $\gamma_{n,i}^{-1}\xo_i$ converge to points in $\rand_i$ as $n\to\infty$. 
Then given arbitrarily small distinct neighborhoods $U_i^+, U_i^-\subset\ganz_i$ of $h_i^+$, $h_i^-$, $i=1,2,\ldots,r$, there exist a finite set $\Lambda\subset \langle g,h\rangle$ and $N_0\in\NN$ \st for some fixed $\lambda,\mu\in\Lambda$ and $n>N_0$ the isometries $\varphi_n:=\lambda\gamma_n\mu^{-1}$ are all regular axial with 
attractive and repulsive fixed points
$\widetilde{\varphi_n}^{+}$, $\widetilde{\varphi_n}^{-}\in\regrand$ satisfying
\[ 
\pi^F(\widetilde{\varphi_n}^{+})\in U_1^+\times U_2^+\times\cdots\times U_r^+,\quad 
\pi^F(\widetilde{\varphi_n}^{-})\in U_1^-\times U_2^-\times\cdots\times U_r^-.\]
\end{cor}
The next result in this section will be the main tool for the construction of certain free subgroups according to Proposition~\ref{constrfreesubgrps}. 
\begin{prp}\label{constraxial}
Assume that $g\in\is(\XX_1)\times\is(\XX_2)\times\cdots\times\is(\XX_r)$ is regular axial and\break $g$, $h$ are independent. Fix a regular axial isometry $\beta\in\is(\XX_1)\times\is(\XX_2)\times\cdots\times\is(\XX_r)$ 
and let ${\cal C}\subset\RR_{> 0}^r$ be an open cone containing $L(\beta)$.
Then for all neighborhoods $V_i^+$, $V_i^-\subset\ganz_i$ of $g_i^+$, $g_i^-$, $i=1,2,\ldots,r$, there exists a regular axial isometry $\alpha\in \langle g,h,\beta\rangle$ with  
\[L(\alpha)\in{\cal C},
\quad \alpha_i^+\in V_i^+\quad\text{and}\quad \alpha_i^-\in V_i^-\quad\text{for all}\quad i\in\{1,2,\ldots,r\}.\] 
\end{prp}
\prf\  Fix $i\in\{1,2,\ldots,r\}$. For $\eta_i\in \{g_i^-,g_i^{+},h_i^-,h_i^{+}\}$ we let $V_i(\eta_i)\subset U_i(\eta_i)\subset \ganz_i$ be neighborhoods of $\eta_i$ and $c_i>0$ as in Proposition~\ref{combgeomlength}.  Making $V_i(g_i^\pm)$ and $U_i(g_i^\pm)$ smaller if necessary we may further assume that 
\[ V_i(g_i^+)\subset U_i(g_i^+)\subset V_i^+\quad\text{and}\quad V_i(g_i^-)\subset U_i(g_i^-)\subset V_i^-.\]
Since the sequences $\beta^n\xo$ and $(\beta^n)^{-1}\xo=\beta^{-n}\xo$ converge to the attractive and repulsive fixed points $\widetilde{\beta^+}$, $\widetilde{\beta^-}\in\regrand$ of $\beta$, Corollary~\ref{gettingregularaxialsinGamma} provides a finite set $\Lambda\subset \langle g,h\rangle$,  $\lambda,\mu\in \Lambda$ 
and $N_0\in\NN$ \st for  all $n> N_0$ isometries $\varphi_n:=\lambda \beta^n\mu^{-1}$ are regular axial with 
\begin{align*} & 
\pi^F(\widetilde{\varphi_n}^{+})\in V_1(h_1^+)\times V_2(h_2^+)\times\cdots\times V_r(h_r^+),\quad\\
& 
\pi^F(\widetilde{\varphi_n}^{-})\in V_1(h_1^-)\times V_2(h_2^-)\times\cdots\times V_r(h_r^-).
\end{align*}
We set $c:=\max\big\{c_i: i\in\{1,2,\ldots,r\}\big\}$, $b:=\max\big\{ d_i(\xo_i,\Ax(\beta_i)):i\in\{1,2,\ldots,r\}\big\}$ and
\[d:=\max\{  d_i(\xo_i,\lambda_i\xo_i):i\in\{1,2,\ldots,r\}, \lambda=(\lambda_1,\lambda_2,\ldots,\lambda_r)\in\Lambda\}.\]
Then we get for $n> N_0$, $k\in\NN$ and $i\in\{1,2,\ldots,r\}$ 
\begin{eqnarray*}
 l_i(\varphi_n^k) &\le &d_i(\xo_i,\varphi_{n,i}^k\xo_i)\le  l_i(\varphi_n^k) +2c,\\
 l_i(\beta^n) &\le & d_i(\xo_i,\beta_i^n\xo_i)\le l_i(\beta^n)+2b\nonumber
 \end{eqnarray*}
and
\begin{align*}
|d_i(\xo_i,\varphi_{n,i}\xo_i)-d_i(\xo_i,\beta_i^n\xo_i)|\le  d_i(\xo_i, \lambda_i\xo_i)+d_i(\xo_i,\mu_i^{-1}\xo_i)\le 2 d,
\end{align*}
which gives (in the special case $k=1$)
\begin{equation}\label{varphibeta}
 |l_i(\varphi_n)-l_i(\beta^n)|\le 2b+2c+2d.
 \end{equation}
We now fix $n>N_0$ and write $\varphi_n=(\varphi_{n,1},\varphi_{n,2},\ldots,\varphi_{n,r})$ . Since for $i\in\{1,2,\ldots,r\}$  we have $\varphi_{n,i}^\pm \in V_i(h_i^\pm)$, Proposition~\ref{combgeomlength} (b) implies the existence of $N, N_n\in \NN$ \st the isometry 
\[ \gamma_n:=g^N \varphi_{n}^{N_n} g^{N}\] 
satisfies
\[ | l_i(\gamma_n)-2l_i(g^N)- l_i(\varphi_n^{N_n})|\le 4c\cdot 3=12c\]
for all $i\in\{1,2,\ldots,r\}$.
Using~(\ref{varphibeta}), $l_i(\varphi_n^{N_n})=N_n\cdot l_i(\varphi_n)$ and $l_i(\beta^n)=n\cdot l_i(\beta)$ we get
\[ |  l_i(\gamma_n)-2 l_i(g^N) -n N_n\cdot  l_i(\beta) |\le 2 N_n(b+c+d)+12c,\]
which implies 
\[ \lim_{n\to\infty} \frac{ l_i(\gamma_n)}{n N_n} =l_i(\beta).\]
This shows that for $n$ sufficiently large we have $ L(\gamma_n)\in{\cal C}$. 
 \qed\\[3mm]
For the last result in this section we assume that $\Gamma<\is(\XX_1)\times\is(\XX_2)\times\cdots\times\is(\XX_r)$ is a discrete group which contains a pair of independent regular axial isometries\break $g=(g_1,g_2,\ldots, g_r)$ and $h=(h_1,h_2,\ldots,h_r)$. As before we fix a base point\break $\xo=(\xo_1,\xo_2,\ldots,\xo_r)\in\Ax(h)$. 

We are going to  construct a generic product for $\Gamma$ as in \cite{MR1933790}, Proposition~2.3.1, which is the essential tool in the proof of Theorem~\ref{concave}. The idea behind is to find a finite set in $\Gamma\times \Gamma$ which maps pairs of orbit points $(\alpha\xo, \beta^{-1}\xo)$ close to a set $\Ax(g)$ or $\Ax(h)$. 
Unfortunately, unlike in the case of symmetric spaces, we do not dispose of 
an equivalent of the result of Abels-Margulis-Soifer (\cite[Proposition~2.3.4]{MR1933790}) which plays a crucial role in the article by Quint. Instead, as in Section~7 of \cite{MR2629900} we will exploit  the dynamics of a free subgroup in $\langle g, h\rangle<  
\Gamma$.

\begin{prp}\label{genprod}
There exists a map $\pr:\Gamma\times\Gamma\to \Gamma$ with the following properties:
\begin{enumerate}
\item[(a)] There exists $\kappa\ge 0$ \st for all $\alpha,\beta\in\Gamma$ we have
$$ \Vert H\big(\pr(\alpha,\beta)\big)-H(\alpha)-H(\beta)\Vert\le \kappa\,.$$ 
\item[(b)] For any $t>0$ there exists a finite set $\Lambda\subset \Gamma$ \st for all $\alpha,\beta,\widehat\alpha,\widehat\beta\in \Gamma$ with $\Vert H(\alpha)-H(\widehat\alpha)\Vert\le t$, $\Vert H(\beta)-H(\widehat\beta)\Vert\le t$ we have
$$ \pr(\alpha,\beta)=\pr(\widehat\alpha,\widehat\beta)\quad\Longrightarrow\qquad \widehat\alpha\in \alpha \Lambda \ \an \ \widehat\beta\in \Lambda \beta\,.$$
\end{enumerate}
\end{prp}
\prf \ 
%
In order to construct a map satisfying property~(a) we  let $\alpha=(\alpha_1,\alpha_2,\ldots,\alpha_r)$, $\beta=(\beta_1,\beta_2,\ldots,\beta_r)\in\Gamma$ arbitrary. 

For $i\in\{1,2,\ldots,r\}$ and $\eta_i\in \{g_i^-,g_i^{+},h_i^-,h_i^{+}\}$ we let $U_i(\eta_i)\subset \ganz_i$ be the neighborhoods of $\eta_i$ and $c_i>0$ the constant provided by Proposition~\ref{combgeomlength}. According to Proposition~\ref{movingpointsnbhds} there exist a finite set $\Lambda\subset\Gamma$ and $\mu=\mu(\alpha)$, $\lambda=\lambda(\beta) \in\Lambda$ \st
\begin{align*} 
  \mu\alpha^{-1}\xo & \in U_1(h_1^-)\times U_2(h_2^-)\times\cdots\times U_r(h_r^-)\quad\text{and}\\
 \lambda\beta\xo & \in U_1(h_1^+)\times U_2(h_2^+)\times\cdots\times U_r(h_r^+) .
\end{align*}
We next set $c:=\max\big\{ c_i:i\in\{1,2,\ldots,r\}\big\}$, 
\[d:=\max\{d_i(\xo_i,\lambda_i\xo_i): i\in\{1,2,\ldots,r\}, \lambda=(\lambda_1,\lambda_2,\ldots,\lambda_r)\in\Lambda \}\]
and fix $i\in\{1,2,\ldots,r\}$. 
Since $d_i(\alpha_i\mu_i^{-1}\lambda_i\beta_i\xo_i,\xo_i)= d_i(\lambda_i\beta_i\xo_i,\mu_i\alpha_i^{-1}\xo_i)$,  Proposition~\ref{combgeomlength} (a) implies
\[ | d_i(\alpha_i\mu_i^{-1}\lambda_i\beta_i\xo_i,\xo_i)- d_i(\lambda_i\beta_i\xo_i,\xo_i)-d_i(\xo_i, \mu_i\alpha_i^{-1}\xo_i)|\le 2c;\]
using $d_i(\lambda_i\beta_i\xo_i,\xo_i)=d_i(\beta_i\xo_i,\lambda_i^{-1}\xo_i)$ and
$d_i(\xo_i, \mu_i\alpha_i^{-1}\xo_i)=d_i(\mu_i^{-1}\xo_i,\alpha_i^{-1}\xo_i)$ we
conclude
\[ | d_i(\alpha_i\mu_i^{-1}\lambda_i\beta_i\xo_i,\xo_i)- d_i(\beta_i\xo_i,\xo_i)-d_i(\xo_i, \alpha_i^{-1}\xo_i)|\le 2c+2d.\]
This implies 
\[  \Vert H(\alpha\mu^{-1}\lambda\beta)-H(\beta)-H(\alpha)\Vert  \le \sqrt{r} (2c+2 d)=:\kappa,\]
hence the assignment $\pr(\alpha,\beta):=\alpha\mu(\alpha)^{-1}\lambda(\beta)\beta$ satisfies  property (a).

It remains to prove that the map $\pr$ from above also satisfies property (b). Suppose there exists $t>0$ \st  for any finite set $\Lambda_n\subset\Gamma$  there exist $\alpha_n,\beta_n,\widehat\alpha_n,\widehat\beta_n\in\Gamma$ with 
\begin{align*}
&\Vert H(\alpha_n)-H(\widehat\alpha_n)\Vert\le t,\quad \Vert H(\beta_n)-H(\widehat\beta_n)\Vert\le t\quad\text{and}\\
\tag{$*$} & \pr(\alpha_n,\beta_n)=\pr(\widehat\alpha_n,\widehat\beta_n), \quad\text{but}\quad \alpha_n^{-1}\widehat\alpha_n\notin \Lambda_n\quad\text{or}\quad \widehat\beta_n\beta_n^{-1}\notin \Lambda_n.
\end{align*}  
For $n\in\NN$ we will work here with the finite set 
\[\Lambda_n:=\{\gamma\in\Gamma:\; d(\xo,\gamma\xo)\le n\}\] 
and fix $\alpha_n=(\alpha_{n,1},\alpha_{n,2},\ldots,\alpha_{n,r})$, $\widehat\alpha_n=(\widehat\alpha_{n,1},\widehat\alpha_{n,2},\ldots,\widehat\alpha_{n,r})$,  $\beta_n=(\beta_{n,1},\beta_{n,2},\ldots,\beta_{n,r})$,  $\widehat\beta_n=(\widehat\beta_{n,1},\widehat\beta_{n,2},\ldots,\widehat\beta_{n,r})\in\Gamma$ \st $(*)$ is satisfied. 

Passing to subsequences if necessary we may assume that for all $i\in\{1,2,\ldots,r\}$ the sequences $(\alpha_{n,i}^{-1}\xo_i)$, $(\widehat\alpha_{n,i}^{-1}\xo_i)$, $(\beta_{n,i}\xo_i)$, $(\widehat\beta_{n,i}\xo_i)\subset\XX_i$ converge. Notice that the limit can be a point in $\XX_i$ or in the geometric boundary $\rand_i$. In any case 
Proposition~\ref{movingpointsnbhds} shows that there exist a finite set $\Lambda\subset\Gamma$ and
$\mu$, $\widehat\mu$, $\lambda$, $\widehat \lambda \in\Lambda$ \st for all\break $i\in\{1,2,\ldots,r\}$ and $n\in\NN$ sufficiently large 
\begin{equation}\label{elementinmenge}
\mu_i\alpha_{n,i}^{-1}\xo_i,\, \widehat\mu_i\widehat\alpha_{n,i}^{-1}\xo_i \in U_i(h_i^-) \quad \text{and}\quad 
\lambda_i\beta_{n,i}\xo_i, \, \widehat\lambda_i\widehat\beta_{n,i}\xo_i\in U_i(h_i^+).
\end{equation}
For $n\in\NN$ and $i=1,2,\ldots,r$ we denote $x_{n,i}$ a point on the geodesic path from $\mu_i\alpha_{n,i}^{-1}\xo_i$ to $\lambda_i\beta_{n,i}\xo_i$, and $\widehat{x}_{n,i}$ a point on the geodesic path from  $\widehat\mu_i\widehat\alpha_{n,i}^{-1}\xo_i$ to $\widehat\lambda_i\widehat\beta_{n,i}\xo_i$ \st $d_i(\xo_i,x_{n,i})\le c$ and $d_i(\xo_i,\widehat x_{n,i})\le c$. Furthermore, setting  
 $\gamma_n:=\alpha_n\mu^{-1}\lambda\beta_n=\widehat\alpha_n\widehat\mu^{-1}\widehat\lambda \widehat\beta_n$ and  denoting  $\sigma_{n,i}$, $i=1,2,\ldots,r$,   the geodesic path $\sigma_{\xo_i,\gamma_{n,i}\xo_i}$ there exist $t_i, \widehat t_i>0$ \st 
 \be d_i(\alpha_{n,i}\mu^{-1}_i\xo_i, \sigma_{n,i}(t_i))&=& d_i(\alpha_{n,i}\mu_i^{-1}\xo_i, \sigma_{n,i})=d_i(\xo_i, \mu_i\alpha_{n,i}^{-1}\sigma_{n,i})\le d_i(\xo_i,x_{n,i})\le c,\\
  d_i(\widehat \alpha_{n,i}\widehat\mu^{-1}_i\xo_i, \sigma_{n,i}(\widehat t_i))&=& d_i(\widehat\alpha_{n,i}\widehat\mu_i^{-1}\xo_i, \sigma_{n,i})=d_i(\xo_i, \widehat\mu_i\widehat\alpha_{n,i}^{-1}\sigma_{n,i})\le d_i(\xo_i,\widehat x_{n,i})\le c
 \ee
by~(\ref{elementinmenge}) and Proposition~\ref{combgeomlength} (a). Hence using again 
\[d=\max\{d_i(\xo_i,\lambda_i\xo_i): i\in\{1,2,\ldots,r\}, \lambda=(\lambda_1,\lambda_2,\ldots,\lambda_r)\in\Lambda \}\]
we estimate
\be d_i(\alpha_{n,i}\xo_i,\sigma_{n,i})&\le & d_i(\alpha_{n,i}\xo_i, \alpha_{n,i}\mu_i^{-1}\xo_i) + d_i(\alpha_{n,i}\mu_i^{-1}\xo_i, \sigma_{n,i}(t_i))\le d+c\,,\\
d_i(\widehat \alpha_{n,i}\xo_i,\sigma_{n,i})&\le & d_i(\widehat \alpha_{n,i}\xo_i, \widehat\alpha_{n,i}\widehat\mu_i^{-1}\xo_i) + d_i(\widehat \alpha_{n,i}\widehat\mu_i^{-1}\xo_i, \sigma_{n,i}(\widehat t_i))\le d+c.
\ee
For $n\in\NN$ and $i=1,2,\ldots,r$  we let $y_{n,i}, \widehat y_{n,i}\in\XX_i$ be the points on the geodesic path $\sigma_{n,i}$ \st $d_i(\xo_i, y_{n,i})=d_i(\xo_i,\alpha_{n,i}\xo_i)$ and $d_i(\xo_i, \widehat y_{n,i})=d_i(\xo_i,\widehat\alpha_{n,i}\xo_i)$. Since\break
 $\Vert H( \alpha_n)-H(\widehat\alpha_n)\Vert\le t$ we have $d_i(y_{n,i},\widehat y_{n,i})\le t$, and, by elementary geometric estimates, 
$$d_i(\alpha_{n,i}\xo_i,y_{n,i})\le 2(d+c)\quad\an\qquad d_i(\widehat\alpha_{n,i}\xo_i,\widehat y_{n,i})\le 2(d+c).$$ 
We summarize 
\[ d_i(\xo_i,\alpha_{n,i}^{-1}\widehat \alpha_{n,i}\xo_i) \le  d_i(\alpha_{n,i}\xo_i,y_{n,i})+d(y_{n,i},\widehat y_{n,i})+d_i(\widehat y_{n,i},\widehat\alpha_{n,i}\xo_i)\le 4(d+c)+t,\]
hence $\hspace*{3cm} d(\xo,\alpha_n^{-1}\widehat \alpha_n\xo)\le  \sqrt{r}(4d+4c+t)=:R.$

In particular, for $n>R$ we have $\alpha_n^{-1}\widehat\alpha_n\in \Lambda_n$, and, in order to obtain the desired contradiction, it remains to prove that $\widehat\beta_n\beta_n^{-1}\in \Lambda_n$ for $n$ sufficiently large. 

Notice that $\widehat\beta_n=\widehat\lambda^{-1}\widehat\mu\widehat\alpha_n^{-1}\gamma_n=\widehat\lambda^{-1}\widehat\mu\widehat\alpha_n^{-1}\alpha_n\mu^{-1}\lambda\beta_n$, hence
\begin{eqnarray*}
d(\xo,\widehat\beta_n\beta_n^{-1}\xo)&=& d(\xo, \widehat\lambda^{-1}\widehat\mu\widehat\alpha_n^{-1}\alpha_n\mu^{-1}\lambda\xo)\le \overbrace{d(\xo,\widehat\lambda^{-1}\xo)}^{\le\sqrt{r}d}+\overbrace{d(\widehat\lambda^{-1}\xo,\widehat\lambda^{-1}\widehat\mu\xo)}^{\le\sqrt{r}d}+\\
&&d(\widehat\lambda^{-1}\widehat\mu\xo,\widehat\lambda^{-1}\widehat\mu\widehat\alpha_n^{-1}\alpha_n\xo)+ d(\widehat\alpha_n^{-1}\alpha_n\xo,\widehat\alpha_n^{-1}\alpha_n\mu^{-1}\xo)+d(\mu^{-1}\xo,\mu^{-1}\lambda\xo)\\
&\le&d(\xo,\alpha_n^{-1}\widehat\alpha_n\xo) + 4\sqrt{r}d\le R+ 4\sqrt{r}d. \end{eqnarray*}
This finishes the proof. \qed

\section{The structure of the limit set}\label{limsetdiscreteisom}

The {\hl geometric limit set} of a discrete group $\Gamma$ acting 
by isometries on a locally compact 
Hadamard space is defined by $\Lim:=\overline{\Gamma\at x}\cap\rand$, where $x\in \XX$ is arbitrary. In this section we are going to describe the structure of the geometric limit set for certain groups $\Gamma<\is(\XX_1)\times\is(\XX_2)\times\cdots\times\is(\XX_r)<\is(\XX)$ acting properly discontinuously on the product $\XX$ of $r$ locally compact 
Hadamard spaces $\XX_1$, $\XX_2,\ldots,\XX_r$.  For convenience the {\hl Furstenberg limit set} of $\Gamma$ is defined by $\Flim:=\pi^F(\Lim\cap\regrand)$.  
Moreover, we let
\[P_\Gamma:=\{ \theta\in E: \Lim\cap\rand_\theta\ne \emptyset\} \subset E\]
be the set of all slopes of geometric limit points, and 
$P^{reg}_\Gamma=P_\Gamma\cap E^+$  the set of slopes of regular limit points.

 In \cite{MR2629900} -- when dealing with only two factors -- we were able to prove Theorems~\ref{minimal} and \ref{Product} in the more general context of discrete isometry groups containing a regular axial isometry with projections which do not globally fix a point in the geometric boundary of the corresponding factor and which possess infinitely many limit points.  Unfortunately, the methods used there and in particular Lemma~4.1 of \cite{MR2629900} are not available in the setting of more factors under the above weak assumption.

From here on we therefore assume -- as in the second part of the aforementioned article -- that $\Gamma<\is(\XX_1)\times\is(\XX_2)\times\cdots\times\is(\XX_r)$ {\hl acts properly discontinuously} on $\XX$ and {\hl possesses two independent regular axial isometries}. This requires in particular that all factors of $\XX$ are rank one spaces as for example 
universal covers of geometric rank one manifolds and CAT$(-1)$-spaces such as locally finite trees or manifolds of pinched negative curvature. Moreover -- as already mentioned in the introduction -- every finite-dimensional unbounded locally compact
CAT(0)-cube complex with an essential and cocompact action of its automorphism group can be decomposed into irreducible factors which are either rank one or Euclidean (compare also \cite[Corollary 2.6]{1105.1675}); hence such CAT$(0)$-cube complexes without Euclidean factors constitute interesting examples for our setting.  
As in the previous section we let $h=(h_1,h_2,\ldots, h_r)$ and $g=(g_1,g_2,\ldots, g_r)\in\Gamma$ be independent regular axial elements of $\is(\XX_1)\times\is(\XX_2)\times\cdots\times \is(\XX_r)$ and fix a base point $\xo=(\xo_1,\xo_2,\ldots,\xo_r)\in\Ax(h)$. Recall that $\widetilde{h^+}, \widetilde{h^-},\widetilde{h^+}, \widetilde{h^-}\in\regrand$ are the attractive and repulsive fixed points, and $h^+,h^-,g^+,g^-\in\Frand$ their images by the Furstenberg projection $\pi^F$.

The following important theorem implies that $\Flim$ can be covered by
finitely many $\Gamma$-trans\-lates of an appropriate open set in $\Frand$. 
\begin{thr}\label{minimal}
The Furstenberg limit set is minimal, i.e. $\Flim$ is the smallest non-empty, $\,\Gamma$-invariant closed subset of $\Frand$.
\end{thr}
\prf\ We first show that every non-empty, $\Gamma$-invariant closed subset $A\subset\Frand$ contains $h^+=(h_1^+,h_2^+,\ldots,h_r^+)$: 
Indeed, if $\xi=(\xi_1,\xi_2,\ldots,\xi_r)\in A$ is arbitrary, then according to Proposition~\ref{movingpointsnbhds} there exists $\lambda\in\Gamma$ \st $\lambda\xi\in\horF(h^-)$. 
So $h^n\lambda\xi$ converges to 
$h^+$ as $n\to\infty$ and -- since $A$ is $\Gamma$-invariant and closed -- 
$h^+$ belongs to $A$. 

It remains to prove that  $\Flim=\overline{\Gamma\cdot h^+}$,
so we let $\eta=(\eta_1,\eta_2,\ldots,\eta_r)\in\Flim$ arbitrary.  Since $\eta\in\Flim$, there exists a sequence $(\gamma_n)\subset\Gamma$ \st $\gamma_n\xo$ converges to a point $\tilde\eta\in\Lim\cap\regrand$ with $\pi^F(\tilde\eta)=\eta$. Passing to a subsequence if necessary, we may assume that $\gamma_n^{-1}\xo$ converges to a point $\tilde\zeta\in\Lim\cap\regrand$ and we set $\zeta:=\pi^F(\tilde\zeta)\in\Flim$.  

We first treat the case $\zeta\in\horF(h^+)$. Let $T\gg1$ and  $\eps>0$ be arbitrary. Then Lemma~\ref{convergence} implies the existence of $N\in\NN$ \st for $i\in\{1,2,\ldots,r\}$ and all $n\ge N$ and $t\in[0,T]$ we have 
\begin{align*} d_i(\sigma_{\xo_i,\gamma_{n,i}\xo_i}(t),\sigma_{\xo_i,\gamma_{n,i} h_i^+}(t))=d_i(\sigma_{\gamma_{n,i}^{-1}\xo_i,\xo_i}(t),\sigma_{\gamma_{n,i}^{-1}\xo_i,h_i^+}(t))&\le \frac{\eps}2.
\end{align*}
Moreover, according to Lemma~\ref{prodtopologyreg} the sequences $\gamma_{n,i}\xo_i$ converges to $\eta_i$ for all $i\in\{1,2,\ldots,r\}$, so we also have
\[
d_i(\sigma_{\xo_i,\gamma_{n,i}\xo_i}(t),\sigma_{\xo_i,\eta_i}(t)) \le\frac\eps2
\]
for $t\in [0,T]$ and $n$ sufficiently large.
Hence we conclude that as $n\to\infty$ we have $\gamma_{n,i}h_i^+\to\eta_i$, and therefore 
$\eta\in\overline{\Gamma\cdot h^+}$.

It remains to deal with the case $\zeta\notin\horF(h^+)$.  
Applying Proposition~\ref{movingpointsnbhds} with $h$ replaced by $h^{-1}$ there  exists $\mu\in\Gamma$ \st $\mu\zeta\in \horF(h^+).$
Since $\gamma_n \mu^{-1}\xo$ still converges to $\tilde\eta$, and 
\[ (\gamma_n \mu^{-1})^{-1}\xo=\mu\gamma_n^{-1}\xo\to \mu\zeta\in\horF(h^+)\quad\text{as}\quad n\to\infty,\] after replacing the sequence $\gamma_n$ by $\gamma_n \mu^{-1}$ we are in the first case. \qed 

\begin{thr}\label{Product}
 The regular geometric limit set $\Lim\cap\regrand$ is isomorphic to the product $\Flim\times P_\Gamma^{reg}$. 
\end{thr}
\prf\  If $\tilde\xi\in \Lim\cap\regrand$, then $\pi^F(\tilde\xi)\in F_\Gamma$, and by definition of $P_\Gamma^{reg}$ the slope of $\tilde\xi$ belongs to $P_\Gamma^{reg}$. 

Conversely, let $\eta=(\eta_1,\eta_2,\ldots,\eta_r)\in F_\Gamma$ and $\theta\in P_\Gamma^{reg}$. We have to show that there exists a limit point $\tilde\eta\in\Lim\cap\rand_\theta$ of slope $\theta$ with $\pi^F(\tilde\eta)=(\eta_1,\eta_2,\ldots,\eta_r)$. By definition
of $P_\Gamma^{reg}$ and Lemma~\ref{prodtopologytheta} there
exists a sequence $(\gamma_n)\subset\Gamma$  \st $\theta^{(n)}:=\thet(\xo,\gamma_n\xo)$ converges to $\theta$ as $n\to\infty$. 
Moreover, by compactness of $\Frand=\rand_1\times\rand_2\times\cdots\times\rand_r$ a subsequence of $\gamma_n\xo$ converges to a point $\tilde \xi\in\Lim\cap\rand_\theta$. We set $\xi:=\pi^F(\tilde\xi)$ and notice that $\tilde\eta\in\regrand$ is  the unique point in $(\pi^F)^{-1}(\eta)$ of slope $\theta$.

By {\rm Theorem}~\ref{minimal}  $\,F_\Gamma=\overline{\Gamma \cdot \xi}\,$ is
minimal and closed under the action of $\Gamma$, hence $\overline{\Gamma \cdot \xi}=F_\Gamma$ and therefore
$$\eta\in\overline{\Gamma\cdot\xi}=\pi^F(\overline{\Gamma\at\tilde\xi}).$$
Since the action of $\Gamma$ on the geometric boundary does not change
 the slope of a point, we conclude 
that  the closure of $\Gamma\cdot \tilde\xi$ contains $\tilde\eta$. In particular
$\tilde\eta \in\overline{\Gamma\at\tilde\xi}\subset\Lim.\ $\qed\\[3mm]
\begin{rmke}
Theorems~\ref{minimal} and \ref{Product} remain true under the weaker assumption that 
$\Gamma$ contains a 
regular axial isometry $h$  and that for any $\xi\in\Flim$ and $\eta\in \{h^+,h^-\}$ there exists $
\lambda\in\Gamma$ \st $\lambda\xi\in\horF(\eta)$. 

When only two factors are present, Lemma~4.1 in \cite{MR2629900} shows that this condition is satisfied if $\Gamma$  contains a regular axial isometry and if the projections of $\Gamma$ to $\is(\XX_1)$ and $\is(\XX_2)$ 
do not globally fix a point in $\rand_1$, $\rand_2$  
and possess infinitely many limit points. 
\end{rmke}\\

In the sequel we will establish an important property of the action of $\Gamma$ on the whole Furstenberg boundary $\Frand$,  namely the fact that the Furstenberg limit set $\Flim$ is a so-called {\hl boundary limit set} for the action of $\Gamma$ on $\Frand$ (see \cite[Definition~4.2]{1105.1675} and also \cite{MR2087783}): (1) clearly implies minimality of $\Flim$ and therefore Theorem~\ref{minimal}; with the notions from \cite[Chapter VI, (1.2)]{MR1090825} (2) says that $\Frand$ and the action of $\Gamma$ on $\Frand$ are {\hl proximal}, and  (3) states that every open set in $\Frand$ is {\hl contractible} to a point in $\Flim$ .
\begin{thr}\label{boundlimset}
The $\,\Gamma$-invariant subset $\Flim\subset\Frand$ satisfies the following:
\begin{itemize}
\item[(1)] 
For all $\zeta\in\Frand$ and every open subset $U\subset\Frand$ with $U\cap \Flim\ne\emptyset$ there exists $\gamma\in\Gamma$ \st $\gamma\zeta\in U$. 
\item[(2)] 
For all $\zeta,\eta\in\Frand$ there exists $\xi\in\Flim$ \st for any neighborhood $U$ of $\xi$ there exists $\gamma\in\Gamma$ with $\{\gamma \zeta,\gamma\eta\}\subset U$.
\item[(3)] 
For all $\zeta\in\Frand$ there exists a neighborhood $V$ of $\zeta$ and a point $\xi\in\Flim$ \st for any neighborhood $U$ of $\xi$ there exists $\gamma\in\Gamma$ with $\gamma V\subset U$.
\end{itemize}
\end{thr}
\prf\ In order to prove (1) we let $\zeta=(\zeta_1,\zeta_2,\ldots, \zeta_r)\in\Frand$ arbitrary and choose an open subset $U\subset\Frand$ with $U\cap\Flim\ne\emptyset$. Let $\xi=(\xi_1,\xi_2,\ldots,\xi_r)\in U\cap\Flim$, and $U_i\subset\rand_i$ open neighborhoods of $\xi_i$, $i=1,2,\ldots,r$, \st \mbox{$U_1\times U_2\times\cdots\times U_r\subset U$.}  Since $\xi\in\Flim$, there exists a sequence $(\gamma_n)\subset\Gamma$ \st $\gamma_n\xo$ converges to a point $\tilde\xi\in\regrand$ with $\pi^F(\tilde\xi)=\xi$. Moreover, passing to a subsequence if necessary we can assume that $\gamma_n^{-1}\xo$ converges to a point $\tilde\eta\in\regrand$, and we set $\eta=\pi^F(\tilde\eta)$. 
If\break $\eta\notin\horF(h^+)$, Proposition~\ref{movingpointsnbhds} with $h$ replaced by $h^{-1}$ provides an element $\mu\in\Gamma$ \st $\mu\cdot\eta\in \horF(h^+).$
Replacing the sequence $\gamma_n$ by $\gamma_n \mu^{-1}$ if necessary we can therefore assume that $\eta\in\horF(h^+)$. 
Now we conclude as in the proof of Theorem~\ref{minimal} 
that for all $i\in\{1,2,\ldots,r\}$  the sequence 
$\gamma_{n,i}h_i^+$ converges to $\xi_i$ as $n\to\infty$; in particular, for some fixed and sufficiently large $n\in\NN$ the regular axial isometry $\varphi=(\varphi_1,\varphi_2,\ldots,\varphi_n):=\gamma_n h\gamma_n^{-1}$ satisfies $\varphi^+\in U_1\times U_2\times\cdots\times U_r\subset U$. According to Lemma~\ref{separatepoints} there exists $\alpha\in\Gamma$ \st $\alpha\zeta\in \horF(\varphi^-)$. So
 $\varphi^n\alpha\zeta$ converges to $\varphi^+$ and hence belongs to $U$ for all sufficiently large $n$.

Proposition~\ref{contract2points} shows that (2) holds with $\xi=h^+$ for all $\zeta,\eta\in\Frand$; (3) follows from Proposition~\ref{movingpointsnbhds} (again with $\xi=h^+$). \qed\\[3mm]
The following theorem  can be viewed as a  strong topological version of the double ergodicity property of Poisson boundaries due to Burger-Monod (\cite{MR1911660}) and Kaimanovich (\cite{MR2006560}). For its proof we will need an important definition as a substitute for  the more familiar notion of $\Gamma$-duality used e.g. in \cite{MR1383216} and \cite{MR2585575} when dealing with only one rank one Hadamard space. 
\begin{df}
Two points $\xi=(\xi_1,\xi_2,\ldots,\xi_r)$, $\eta=(\eta_1,\eta_2,\ldots,\eta_r)\in \Frand$ are called {\hd $\,\Gamma$-related} if for all $i\in\{1,2,\ldots,r\}$ and all neighborhoods $U_i$ of $\xi_i$ and  $V_i$ of $\eta_i$  in $\ganz_i$ there exists $\gamma=(\gamma_1,\gamma_2,\ldots,\gamma_r)\in\Gamma$ 
\st  
\[ \gamma_i(\ganz_i\setminus U_i)\subset V_i\quad\text{and}\quad \gamma_i^{-1}(\ganz_i\setminus V_i)\subset U_i\quad\text{for all}\quad i\in\{1,2,\ldots,r\}.\]
We will denote $\rel(\xi)$ the set of points in $\Frand$ which are $\Gamma$-related to $\xi$.
\end{df}
Notice that for any $\xi\in\Frand$ the set $\rel(\xi)$ is closed with respect to the topology of $\Frand$. Moreover, if $\eta\in\rel(\xi)$, then $\eta_i$ is $\Gamma_i$-dual to $\xi_i$ for all $i\in\{1,2,\ldots,r\}$.  The converse clearly does not hold in general.

The importance of the notion lies in the following. If $\widetilde{h^+}$, $\widetilde{h^-}$ denote the attractive and repulsive fixed point of a regular axial isometry 
$h\in\Gamma$, then by  Lemma~\ref{dynrankone} (c) the points $h^+= \pi^F(\widetilde{h^+})$ and $h^-= \pi^F(\widetilde{h^-})$ are $\Gamma$-related. Conversely, if two points $\xi=(\xi_1,\xi_2,\ldots,\xi_r)$, $\eta=(\eta_1,\eta_2,\ldots,\eta_r)\in\Frand$ are $\Gamma$-related, then by definition there exists a sequence $(\gamma_n)=\big((\gamma_{n,1},\gamma_{n,2},\ldots,\gamma_{n,r})\big)\subset\Gamma$ \st for all $i\in\{1,2,\ldots,r\}$ we have $\gamma_{n,i}\xo_i\to\eta_i$ and $\gamma_{n,i}^{-1}\xo_i\to\xi_i$ as $n\to\infty$. Hence if for $i\in\{1,2,\ldots,r\}$ the points $\xi_i$, $\eta_i\in\rand_i$ can be joined by a rank one geodesic, then in view of  Lemma~\ref{elementsarerankone} $\gamma_{n,i}$ is rank one for $n$ sufficiently large and satisfies 
$$\gamma_{n,i}^+\to\eta_i\quad\mbox{and}\qquad \gamma_{n,i}^{-}\to\xi_i\quad\text{as}\quad  n\to\infty.$$
For the sequel we denote $\Delta \subset\Frand\times\Frand$ the generalized diagonal 
$$\Delta:=\big\{(\xi,\eta)\in \Frand\times\Frand : \xi_i= \eta_i\  \text{ for some }\  i\in\{1,2,\ldots,r\}\big\}\,.$$
With this notation we have the following 
\begin{thr}\label{densityofaxials}
The set of pairs of fixed points $(\gamma^+,\gamma^-)\subset\Frand\times \Frand$ of regular axial isometries  $\gamma\in\Gamma$ is dense in $\big(F_\Gamma
\times F_\Gamma\big)\setminus \Delta$.
\end{thr}
\prf\ Recall that $g=(g_1,g_2,\ldots,g_r)$, $h=( h_1, h_2,\ldots,h_r)\in\Gamma$ are two independent regular axial isometries. In view of the paragraph preceding the theorem we 
first prove that any two distinct points  in $\{g^-,g^+,h^-, h^+\}$ are $\Gamma$-related. 

For $i\in\{1,2,\ldots,r\}$ and $\eta_i\in \{g_i^-,g_i^{+},h_i^-,h_i^{+}\}$ we  let $U_i(\eta_i)\subset \ganz_i$ be an arbitrary, sufficiently small neighborhood of $\eta_i$ with $\xo_i\notin U_i(\eta_i)$ \st all $U_i(\eta_i)$ are pairwise disjoint. According to Lemma~\ref{dynrankone} (c) there exists  a constant   $N\in\NN$ \st for all $i\in\{1,2,\ldots,r\}$
\begin{equation}\label{factorpingpongg}  
g_i^{\pm N}\big(\ganz_i\setminus U_i(g_i^{\mp})\big)\subset U_i(g_i^{\pm})\qquad\text{and}\quad h_i^{\pm N}\big(\ganz_i\setminus U_i(h_i^{\mp})\big)\subset U_i(h_i^{\pm}).
\end{equation}
Let $\gamma,\varphi \in\{g,g^{-1},h,h^{-1}\}$, $\varphi\ne \gamma$.  Using the fact that either $\varphi=\gamma^{-1}$ or $\gamma_i$, $\varphi_i$ are independent for $i\in\{1,2,\ldots,r\}$  property~(\ref{factorpingpongg}) implies 
\[ \gamma_i^N \varphi_i^{-N}\big(\ganz_i\setminus U_i(\varphi_i^+)\big)\subset U_i(\gamma_i^+)\quad \text{and}\quad (\gamma_i^N\varphi_i^{-N})^{-1}\big(\ganz_i\setminus U_i(\gamma_i^+)\big)\subset U_i(\varphi_i^+)\]
 for $i\in\{1,2,\ldots,r\}$. Hence $\varphi^+\in\rel(\gamma^+)$.
 
 Next we will show that any $\xi=(\xi_1,\xi_2,\ldots,\xi_r)\in\Flim$ with $\xi_i\notin\{g_i^-,g_i^+,h_i^-,h_i^+\}$\break for all $i\in\{1,2,\ldots,r\}$  is $\Gamma$-related to an arbitrary point in $\{g^-,g^+,h^-,h^+\}$.  
For\break $i\in\{1,2,\ldots,r\}$ and $\zeta_i\in \{\xi, g_i^-,g_i^{+},h_i^-,h_i^{+}\}$ we let $U_i(\zeta_i)\subset \ganz_i$ be a sufficiently small neighborhood of $\zeta_i$ with $\xo_i\notin U_i(\zeta_i)$  \st all $U_i(\zeta_i)$ are pairwise disjoint.  By  Lemma~\ref{elementsarerankone} there exist neighborhoods $W_i(\zeta_i)\subset U_i(\zeta_i)$, $\zeta_i\in \{\xi, g_i^-,g_i^{+},h_i^-,h_i^{+}\}$, \st every $\gamma_i\in\Gamma_i$ with $\gamma_i\xo_i\in W_i(\zeta_i)$, $\gamma_i^{-1}\xo_i\in W_i(\eta_i)$, $\eta_i\in\{\xi_i, g_i^-,g_i^{+},h_i^-,h_i^{+}\}\setminus\{\zeta_i\}$,  is rank one with $\gamma_i^+\in U_i(\zeta_i)$ and $\gamma_i^-\in U_i(\eta_i)$. 

Since $\xi\in\Flim$, there exists a sequence $(\gamma_n)=\big((\gamma_{n,1},\gamma_{n,2},\ldots,\gamma_{n,r})\big)\subset\Gamma$ \st $\gamma_{n,i}\xo_i\to \xi_i$ for all $i\in\{1,2,\ldots,r\}$. Upon passing to  a subsequence if necessary we may assume that $\gamma_{n,i}^{-1}\xo_i$ converges to a point  in $\rand_i$ for all $i\in\{1,2,\ldots,r\}$.  By Proposition~\ref{movingpointsnbhds} there exist  a finite set $\Lambda\subset\Gamma$ and $\mu\in\Lambda$ \st for all $n$ sufficiently large we have   
\[\mu\gamma_n^{-1}\xo\in W_1(h_1^-)\times W_2(h_2^-)\times\cdots\times W_r(h_r^-).\]  
Moreover, since for $i\in\{1,2,\ldots, r\}$ and $x_i\in\XX_i$ the sequence  $(\gamma_{n,i} x_i)$ converges to $\xi_i$, we also have 
\[ \gamma_n\mu^{-1}\xo\in W_1(\xi_1)\times W_2(\xi_2)\times\cdots\times W_r(\xi_r)\]
for $n$ sufficiently large. 
By  Lemma~\ref{elementsarerankone} we conclude that for $n$ sufficiently large  the isometry $\gamma_n\mu^{-1}$ is regular axial with 
\[(\gamma_n\mu^{-1})^+\in U_1(\xi_1)\times U_2(\xi_2)\times\cdots\times U_r(\xi_r),\quad (\gamma_n\mu^{-1})^-\in  U_1(h_1^-)\times U_2(h_2^-)\times\cdots\times U_r(h_r^-).\]
 This implies that $\xi\in\rel(h^-)$ and by symmetry 
\begin{equation}\label{relatedtofourpoints}
\xi\in\rel(g^-)\cap \rel(g^+)\cap\rel(h^-)\cap\rel(h^+)\,.
\end{equation}

Next we let $\xi=(\xi_1,\xi_2,\ldots,\xi_r)$ and $\eta=(\eta_1,\eta_2,\ldots,\eta_r)\in\Flim$  \st for  all\break $i\in\{1,2,\ldots,r\}$ we have  $\{ \xi_i, \eta_i\}\cap\{g_i^-,g_i^+,h_i^-,h_i^+\}=\emptyset$ and $\xi_i\ne\eta_i$. As above, for 
$\zeta_i\in \{\xi_i,\eta_i,  h_i^-\}$ we let $U_i(\zeta_i)\subset \ganz_i$ be a small neighborhood of $\zeta_i$ with $\xo_i\notin U_i(\zeta_i)$ \st all $U_i(\zeta_i)$ are pairwise disjoint.  
By the arguments in the previous paragraph there exists a regular axial isometry $\varphi\in\Gamma$ with 
\[ \varphi^+\in U_1(\xi_1)\times U_2(\xi_2)\times\cdots\times U_r(\xi_r)\quad\text{and}\quad \varphi^-\in U_1(h_1^{-})\times U_2(h_2^{-})\times\cdots\times U_r(h_r^-).\]
 In particular, $\varphi_i$ and $g_i$ are independent for $i=1,2,\ldots,r$. Replacing $h$ by $\varphi$ in (\ref{relatedtofourpoints}) we know that $\eta\in \rel(g^-)\cap \rel(g^+)\cap\rel(\varphi^-)\cap\rel(\varphi^+)$, in particular $\eta\in\rel(\varphi^+)$. So using the fact that $\eta_i$ can be joined to $\varphi_i^+$ by a rank one geodesic in $\XX_i$ for all\break $i\in\{1,2,\ldots,r\}$, given small neighborhoods $U_i(\varphi_i^+)\subset U_i(\xi_i)$ 
 there exists $\gamma\in\Gamma$ regular axial  with 
 \[\gamma^+\in U_1(\varphi_1^+)\times U_2(\varphi_2^+)\times\cdots\times U_r(\varphi_r^+)\subset U_1(\xi_1)\times U_2(\xi_2)\times\cdots\times U_r(\xi_r)\] 
 and $\gamma^-\in U_1(\eta_1)\times U_2(\eta_2)\times\cdots\times U_r(\eta_r)$. \qed

\section{The limit cone}\label{sectlimcone}

Given a discrete group $\Gamma<\is(\XX_1)\times\is(\XX_2)\times\cdots\times\is(\XX_r)$, the {\hl limit cone} $\ell_\Gamma$ of $\Gamma$  is defined as the closure of the set of half-lines in $\RR_{\ge 0}^r$ spanned by the set of vectors 
\[ \{ L(\gamma)\in\RR_{>0}^r: \gamma
\in\Gamma\  \ \text{regular axial}\}.\] 

Notice that this definition 
differs from the one given by Y.~Benoist in \cite{MR1437472}, where 
 the translation vectors of {\hl all} elements in $\Gamma$ are considered. However,  Y.~Benoist showed that in the case of reductive groups one can equivalently use only the translation vectors of $\RR$-regular elements in the definition of the limit cone; so our definition can be viewed as an appropriate analogous one. 


As before we assume that $\Gamma$ contains a pair of independent regular axial isometries $g=(g_1,g_2,\ldots, g_r), h=(h_1,h_2,\ldots, h_r)$ and fix a base point \mbox{$\xo=(\xo_1,\xo_2,\ldots,\xo_r)\in\Ax(h)$.} 
The following theorem shows that the limit cone is closely related to the set $P_\Gamma$ introduced at the beginning of Section~\ref{limsetdiscreteisom}.

\begin{prp}\label{propertiesofPGamma}
We have the follwing inclusions:
 \[ \ell_\Gamma\cap E\subset P_\Gamma\qquad\text{and}\qquad P_\Gamma^{reg} \subset \ell_\Gamma\cap E^+.\]  
\end{prp}
\prf\ We first show $\ell_\Gamma\cap E\subset P_\Gamma$: If $(\gamma_n)=\big((\gamma_{n,1},\gamma_{n,2},\ldots,\gamma_{n,r})\big)$ is a sequence of regular axial isometries \st $\widehat L(\gamma_n)$  converges to $\theta=(\theta_1,\theta_2,\ldots,\theta_r)\in E$, we choose 
$$k_n\ge 2n \max\Big\{\frac{d(\xo,\Ax(\gamma_{n}))}{l(\gamma_{n})}, \frac{d_i(\xo_i,\Ax(\gamma_{n,i}))}{l_i(\gamma_{n})}:i\in\{1,2,\ldots,r\}\Big\}$$ and set
$\varphi_n:=\gamma_n^{k_n}$. From
$$ k_n l_i(\gamma_{n})\le d_i(\xo_i,\varphi_{n,i}\xo_i)\le k_nl_i(\gamma_{n})+2 d_i(\xo_i,\Ax(\gamma_{n,i})) \le k_n l_i(\gamma_{n})\big(1+\frac1{n}\big)$$
for $i\in\{1,2,\ldots,r\}$ and 
\[  k_n l (\gamma_{n})\le d(\xo,\varphi_{n}\xo)\le k_nl(\gamma_{n}) +2 d(\xo,\Ax(\gamma_{n}))\le k_n l(\gamma_{n})\big(1+\frac1{n}\big),\]
we conclude that  (by definition of $\widehat L(\gamma_n)$)
\begin{align*}
\theta_i & = \lim_{n\to\infty}\left(\frac{l_i(\gamma_{n})}{l(\gamma_n)}\right)= \lim_{n\to\infty}\left(\frac{l_i(\gamma_{n})\cdot \big(1+\frac1{n}\big)}{l(\gamma_{n})}\right) \ge \lim_{n\to\infty}\frac{d_i(\xo_i,\varphi_{n,i}\xo_i)}{d(\xo,\varphi_{n}\xo)},\\
\theta_i & = \lim_{n\to\infty}\left(\frac{l_i(\gamma_{n})}{l(\gamma_{n})\cdot \big(1+\frac{1}{n}\big)}\right) \le \lim_{n\to\infty}\frac{d_i(\xo_i,\varphi_{n,i}\xo_i)}{d(\xo,\varphi_{n}\xo)},
\qquad i=1,2,\ldots, r,
\end{align*}
which shows that $\thet(\xo,\varphi_n\xo)$ converges to $\theta$ as $n\to\infty$. 


Next we prove the inclusion $P_\Gamma^{reg}\subset \ell_\Gamma\cap E^+$: If $\theta\in P_\Gamma^{reg}$, then by definition of $P_\Gamma^{reg}$ there exists a point $\tilde\xi\in\Lim\cap\rand_\theta\subset\regrand$; in particular,    there exists a sequence $(\gamma_n)=\big((\gamma_{n,1},\gamma_{n_2},\ldots,\gamma_{n,r})\big)\subset\Gamma$ \st $\gamma_n\xo$ converges to $\tilde\xi$ 
and hence necessarily the sequence of directions $\thet(\xo,\gamma_{n}\xo)$ converges to $\theta=(\theta_1,\theta_2,\ldots,\theta_r)$  as $n\to\infty$. Passing to a subsequence if necessary we can assume that $\gamma_{n}^{-1}\xo$ converges to a point $\tilde\zeta\in\rand$ (which necessarily also belongs to $\rand_\theta\subset\regrand$)  as $n\to\infty$.  For $i\in\{1,2,\ldots,r\}$ we  let
$V_i(h_i^\pm)\subset\ganz_i$ be neighborhoods of $h_i^\pm$ and $c_i>0$ such that the assertion of Proposition~\ref{combgeomlength} 
holds in the factor $\XX_i$.  By 
Corollary~\ref{gettingregularaxialsinGamma}
there exist  a finite set $\Lambda\subset\Gamma$, $\lambda$, $\mu\in\Lambda$ and $N_0\in\NN$ \st for all $n>N_0$ the isometries
\[\varphi_n:=\lambda\gamma_n\mu^{-1}\] 
are regular axial with 
\[\varphi_n^\pm \in V_1(h_1^\pm)\times V_2(h_2^\pm)\times\cdots\times V_r(h_r^\pm).\]
Put $c:=\max\big\{c_i: i\in\{1,2,\ldots,r\}\big\}$,  
\[d:=\max\big\{d_i(\xo_i,\lambda_i\xo_i): \lambda=(\lambda_1,\lambda_2,\ldots,\lambda_r)\in\Lambda,\  i\in\{1,2,\ldots,r\} \big\}\]
and fix $n>N_0$. Writing $\varphi_n=(\varphi_{n,1},\varphi_{n_2},\ldots,\varphi_{n,r})$ the triangle inequality implies
\begin{equation}\label{samedistance}
|d_i(\xo_i,\varphi_{n,i}\xo_i)-d_i(\xo_i,\gamma_{n,i}\xo_i)|\le 2d
\end{equation}
  for all $i\in\{1,2,\ldots,r\}$; 
clearly this also gives
\[ |d(\xo,\varphi_{n}\xo)-d(\xo,\gamma_{n}\xo)|\le 2d\sqrt{r}.\]
Moreover, by Proposition~\ref{combgeomlength} (a) we have $d_i(\xo_i,\Ax(\varphi_{n,i}))\le c$, hence
\begin{align*} l_i(\varphi_{n}) & \le d_i(\xo_i,\varphi_{n,i}\xo_i)\le l_i(\varphi_{n})+2c\quad\text{for}\quad i\in\{1,2,\ldots,r\}\quad\text{and}\\
 l(\varphi_{n}) &\le d(\xo,\varphi_{n}\xo)\le l(\varphi_{n})+2c\sqrt{r}.
\end{align*}
In combination with~(\ref{samedistance}) we get 
\[\theta_i=\displaystyle\lim_{n\to\infty} \frac{l_i(\varphi_{n})}{l(\varphi_{n})}\quad\text{for}\quad  i\in\{1,2,\ldots,r\},\] and therefore $\quad\theta=\displaystyle\lim_{n\to\infty} \widehat L(\varphi_n)$. \qed\\[3mm]
Notice that Proposition~\ref{propertiesofPGamma} in particular implies
$ P_\Gamma^{reg}=\ell_\Gamma \cap E^+,$ 
which was proved in \cite[Theorem 5.2]{MR2629900} for the special case of only two factors.  We also want to make the following\\[2mm]
\begin{rmke}
Our proof does not give the stronger statement $ P_\Gamma=\ell_\Gamma\cap E\,$ of Theorem~D. 
This is due to the fact that
 a singular limit point can be approached by orbit points $\gamma_n\xo$ for which the projections to one of the factors $\XX_i$ remain at bounded distance of $\xo_i$. 
However, if $\theta\in P_\Gamma \setminus P_\Gamma^{reg}$ and at least {\hl one} point $\tilde\xi\in \rand_\theta\subset\singrand$ is the limit of a sequence $\gamma_n\xo$ for which the  projections to {\hl all} factors leave every bounded set, 
then our proof of the second inclusion together with Corollary~\ref{gettingregularaxialsinGamma}  
shows that $\tilde\xi$ is also the limit of a sequence of regular axial elements of $\Gamma$ and hence $\theta\in\ell_\Gamma$. This is already remarkable because in the original sequence $(\gamma_n)$ the projections could be parabolic or elliptic of infinite order. Using the exponent of growth in Section~7 we will be able to complete the proof of the full statement of Theorem~D.
\end{rmke}\\[3mm]
Our next goal is to describe the limit cone 
more precisely.  The following easy lemma will be useful for proving convexity as stated in Proposition~\ref{limitconeconvex}:
\begin{lem}\label{sectorinell}
If  $\alpha,\beta\in\Gamma$ are independent regular axial elements with $L(\alpha), L(\beta)\in \ell_\Gamma$, then the convex hull of the half-lines determined by $L(\alpha)$ and $L(\beta)$ is contained in $\ell_\Gamma$.
 \end{lem}
 \prf\ 
Since $\alpha, \beta\in\Gamma$ are independent regular axial isometries, Propositon~\ref{combgeomlength} (b) ensures the existence of $c>0$ and $N\in\NN$ \st for all $i\in\{1,2,\ldots, r\}$ and $k,m\in\NN$ we have
\[ |l_i(\alpha^{kN}\beta^{mN})- kN  l_i(\alpha)-m N l_i(\beta)|\le 4c\cdot 2=8c, \]
which immediately implies 
\[\lim_{n\to\infty} \frac{L(\alpha^{k n N}\beta^{m nN})}{nN}= k L(\alpha)+m L(\beta).\]
Since $\ell_\Gamma$ is closed and $L(\alpha^{k n N}\beta^{m nN})\in\ell_\Gamma$ for all $k,m,n\in\NN$ we conclude that for any positive rational number $q\in\QQ$ the half-line determined by $L(\alpha)+q L(\beta)$ belongs to $\ell_\Gamma$.  So the convex hull of the half-lines determined by $L(\alpha)$ and $L(\beta)$ is included in $\ell_\Gamma$, which we wanted to prove.\qed\\[3mm]
The following proposition is the appropriate analogon of Proposition~2.2.7 in \cite{MR1933790} (compare also Proposition~5.1 in \cite{MR1437472}) for our setting.
\begin{prp}\label{constrfreesubgrps}
There is a constant $\kappa > 0$ \st for every open cone \mbox{${\cal C}\subset\RR^r_{>0}$} 
with ${\cal C}\cap \ell_\Gamma\ne\emptyset$ 
there exist independent regular axial isometries $\alpha,\beta\in\Gamma$ with $L(\alpha)$,\break $L(\beta)\in {\cal C}$ \st the semi-group $\langle \alpha,\beta\rangle^+\subset\Gamma$ is free. 

Moreover, if $\Phi:\langle \alpha,\beta\rangle^+\to\RR^r$ is the unique homomorphism of semi-groups sending $\alpha$ to $L(\alpha)$ and $\beta$ to $L(\beta)$, then for any  word $\gamma\in\langle \alpha,\beta\rangle^+$ of length $n\ge 1$ in the generators $\alpha,\beta$ one has
\[ H(\xo,\gamma\xo)\in {\cal C}\quad\text{and}\quad \Vert H(\xo,\gamma\xo)-\Phi(\gamma)\Vert\le \kappa\cdot n.\] 
\end{prp}
\prf\  For $i\in\{1,2,\ldots,r\}$ we fix neighborhoods $V_i(g_i^\pm)$, $V_i(h_i^\pm)\subset\ganz_i$ and $c_i>0$ \st the assertion of Proposition~\ref{combgeomlength} holds in the factor $\XX_i$. We further set\break $c:=\max\big\{c_i:i\in\{1,2,\ldots,r\}\big\}$.

Since ${\cal C} 
\cap\ell_\Gamma\ne\emptyset\, $ 
 there exist $\alpha',\beta'\in\Gamma$ regular axial with $L(\alpha'),L(\beta')\in {\cal C}$. 
Pro\-position~\ref{constraxial} implies that there exist regular axial isometries  $\alpha=(\alpha_1,\alpha_2,\ldots,\alpha_r)$, $\beta=(\beta_1,\beta_2,\ldots,\beta_r)\in\Gamma$ with $L(\alpha),L(\beta)\in {\cal C}$, 
\[ \alpha_i^\pm\in V_i(g_i^\pm)\quad\text{and}\quad \beta_i^\pm\in V_i(h_i^\pm)\quad\text{for all}\quad i\in\{1,2,\ldots,r\}.\]
Obviously every non-trivial linear combination of $L(\alpha)$, 
$L(\beta)\in\RR^r_{>0}$ with non-negative coefficients is included in the sector ${\cal S}\subset \RR^r_{>0}$ spanned by $L(\alpha)$ and $L(\beta)$; in particular, for any $\gamma\in\langle \alpha,\beta\rangle^+\setminus\{\id\}$ we have $\Phi(\gamma)\in{\cal S}$. Since ${\cal C}$ is an open cone containing $L(\alpha)$ and $L(\beta)$, the whole sector ${\cal S}$ is included in ${\cal C}$ and there exists $\eps>0$ \st every unit vector $\widehat H\in\RR^r_{>0}$ with $\Vert \widehat H- \widehat L\Vert <\eps$ for some unit vector $\widehat L\in{\cal S}$ is contained in ${\cal C}$.

Replacing $\alpha$ and $\beta$ by a sufficiently high power if necessary we can assume that 
\begin{equation}\label{minmove}
\min\{ d(\xo,\alpha\xo),d(\xo,\beta\xo)\} > 8c\sqrt{r}\big(1+\frac1{\eps}\big)
\end{equation}
and that 
the assertion of Proposition~\ref{combgeomlength} holds in each factor $\XX_i$, $i\in\{1,2,\ldots,r\}$, with $N_{\alpha_i}=N_{\beta_i}=1$. 
So let $\gamma=s_1^{k_1} s_2^{k_2}\cdots s_n^{k_n}$ be a word in the semi-group $\langle \alpha,\beta\rangle^+$ with 
 $s_j\in \{\alpha,\beta\}$ and $k_j\in\NN$, $j\in\{1,2,\ldots,n\}$; the word length of $\gamma$ then clearly satisfies $\displaystyle\sum_{j=1}^n k_j\ge n$ and we have\vspace{-6mm}
\[ \Phi(\gamma)=\sum_{j=1}^n k_j L(s_j).\]
Proposition~\ref{combgeomlength} (b) further shows that   
\quad $\displaystyle |d_i(\xo_i,\gamma_i\xo_i)-\sum_{j=1}^n k_j l_i(s_j)|\le 4c \cdot n$\\ 
 for all $i\in\{1,2,\ldots,r\}$,  hence
 \begin{align*}
   \Vert H(\gamma)- \Phi(\gamma)\Vert^2 & = \sum_{i=1}^r \Bigl(\big| d_i(\xo_i,\gamma_i\xo_i)-\sum_{j=1}^n k_j l_i(s_j)\big|^2\Bigr)
   \le  r\cdot (4c\cdot n)^2. 
  \end{align*}
So the second assertion is true with $\kappa:= 4c\sqrt{r}$. 

 Concerning the first assertion we remark that the proof of Proposition~\ref{combgeomlength} (b) implies that
 \[ d(\xo,\gamma\xo)\ge \sum_{j=1}^n d(\xo,s_j^{k_j}\xo)-2c\sqrt{r}\cdot n,\]
 hence by~(\ref{minmove})
 \[ d(\xo,\gamma\xo)\ge n\cdot 8c\sqrt{r}\big(1+\frac1{\eps}\big) -2c\sqrt{r}\cdot n>  \frac{8c\sqrt{r}}{\eps}\cdot n=\frac{2\kappa\cdot  n}{\eps}.\]
 So we estimate 
 \begin{align*} 
 \Big\Vert \frac{H(\xo,\gamma\xo)}{d(\xo,\gamma\xo)}-\frac{\Phi(\gamma)}{\Vert\Phi(\gamma)\Vert}\Big\Vert &\le \frac{1}{d(\xo,\gamma\xo)}\Vert H(\xo,\gamma\xo)-\Phi(\gamma)\Vert +\Big\Vert \frac{\Phi(\gamma)}{d(\xo,\gamma\xo)}-\frac{\Phi(\gamma)}{\Vert \Phi(\gamma)\Vert}\Big\Vert\\
 &\le \frac{1}{d(\xo,\gamma\xo)} \cdot \kappa\cdot n + \Vert \Phi(\gamma)\Vert \cdot \Big| \frac1{d(\xo,\gamma\xo)}-\frac1{\Vert \Phi(\gamma)\Vert}\Big|\\
  &\le \frac{\kappa\cdot n}{d(\xo,\gamma\xo)}+
 \frac1{d(\xo,\gamma\xo)}\big|  \Vert \Phi(\gamma)\Vert- d(\xo,\gamma\xo)\big|\le \frac{2\kappa\cdot n}{d(\xo,\gamma\xo)}<\eps,
  \end{align*}
 where  we used the inverse triangle inequality
  \[ \big|  \Vert \Phi(\gamma)\Vert-\Vert H(\xo,\gamma\xo)\Vert \big| \le \Vert \Phi(\gamma)-H(\xo,\gamma\xo)\Vert \le \kappa\cdot n.\]
 So $\Phi(\gamma)\in{\cal S}$ and the choice of $\eps>0$ imply that 
$ H(\xo,\gamma\xo)\in{\cal C}.$\qed


 \begin{prp}\label{limitconeconvex}
The limit cone $\ell_\Gamma$ is convex. 
\end{prp} 
\prf\ Let $L,L'\in \ell_\Gamma$. Using Pro\-position~\ref{constraxial} and arguments as in the proof of the previous proposition there exist two independent regular axial isometries $\alpha,\beta\in\Gamma$ with $L(\alpha)$ and $L(\beta)$ arbitrarily close to the half-lines determined by $L$ and $L'$. 

From Lemma~\ref{sectorinell} we know that the convex hull of the half-lines determined by $L(\alpha)$ and $L(\beta)$ is contained in $\ell_\Gamma$. Since $\ell_\Gamma$ is closed, the same is true for the convex hull of $L$ and $L'$, which  finishes the proof.
 \qed\\

We finally state a result concerning free subgroups of $\Gamma$ which is a corollary of Proposition~\ref{constraxial} and the proof of Proposition~\ref{constrfreesubgrps}:
\begin{prp}\label{constrfreegrps}
For every open cone ${\cal C}\subset\RR^r_{>0}$ 
with ${\cal C}\cap \ell_\Gamma\ne\emptyset$ 
there exists a free subgroup $\,\Gamma'<\Gamma$ \st 
\[\ell_{\Gamma'}\subset {\cal C}.\]
\end{prp}
\prf\  As in the proof of Proposition~\ref{constrfreesubgrps} and due to Proposition~\ref{combgeomlength} (b) there exist independent regular axial isometries $\alpha,\beta\in\Gamma$  with 
 $L(\alpha), L(\beta)\in{\cal C}$, and a constant $c>0$ such that for every cyclically reduced word $\gamma=s_1^{k_1}s_2^{k_2}\cdots s_n^{k_n}$ with $s_j\in\{\alpha,\alpha^{-1},\beta,
 \beta^{-1}\}$ and $k_j\in\NN$, $j\in\{1,2,\ldots, n\}$ we have
 \[ \Vert L(\gamma)-\sum_{j=1}^n k_j L(s_j)\Vert \le 4c\sqrt{r}\cdot n.\]
Since $L(\alpha^{-1})=L(\alpha)$ and $L(\beta^{-1})=L(\beta)$, this shows that $L(\gamma)$ is at distance $\le 4c\sqrt{r}\cdot n$ of a non-trivial linear combination of $L(\alpha)$, 
$L(\beta)\in\RR^r_{>0}$ with non-negative coefficients (which is included in the sector ${\cal S}\subset{\cal C}$ spanned by $L(\alpha)$ and $L(\beta)$). Passing to powers of $\alpha$ and $\beta$ if necessary we can arrange (as in the proof of Proposition~\ref{constrfreesubgrps}) that for every cyclically reduced word  $\gamma\in \langle \alpha,\beta\rangle$ the translation vector $L(\gamma)$ is arbitrarily close to ${\cal S}$ and hence also included in ${\cal C}$.  This finishes the proof, because every element in $\Gamma':=\langle \alpha,\beta\rangle$ is conjugated to a cyclically reduced element as above, and the translation vector is invariant by conjugation.\qed

%

\section{The exponent of growth for a given slope}\label{ExpGrowthSlope}

For the remainder of the article we let $\Gamma<\is(\XX_1)\times\is(\XX_2)\times\cdots\times \is(\XX_r)$ be a group acting properly discontinuously by isometries on a product $\XX=\XX_1\times \XX_2\times\cdots\times \XX_r$  of $r$ locally compact 
Hadamard spaces which contains two independent regular axial isometries $h=(h_1,h_2,\ldots,h_r)$ and $g=(g_1,g_2,\ldots, g_r)$. We also fix a base point\break $\xo=(\xo_1, \xo_2,\ldots,\xo_r)\in\Ax(h) $.

Recall the notation introduced in Section~\ref{prodHadspaces}; in particular, we denote $E\subset \RR^r$ the set of unit vectors in $\RR_{\ge 0}^r$. 
In this section we want to describe the map which assigns to each vector $\theta\in E$  the exponential growth rate $\delta_\theta(\Gamma)$ of orbit points of $\Gamma$ in $\XX$ with  prescribed slope  $\theta$.   For $x,y\in\XX$, $\theta \in E$ and  $\eps>0$ we first set
$$ \Gamma(x,y;\theta,\eps):=\{\gamma\in\Gamma: \gamma y\ne x\quad \mbox{and}\ \ \Vert \thet(x,\gamma y)-\theta \Vert <\eps \}.$$
In order to define $\delta_\theta(\Gamma)$  we first introduce the following partial sum of the Poincar{\'e} series for $\Gamma$. For $s\in\RR$ we consider the series
$$Q^{s,\eps}_{\theta}(x,y)=\sum_{\gamma\in
  \Gamma(x,y;\theta,\eps)}
  \e^{-s d(x,\gamma y)}$$
and denote $\delta_\theta^\eps(x,y)$ its {\hl critical exponent}, i.e. the   unique
  real number \st $Q^{s,\eps}_{\theta}(x,y)$ converges for
  $s>\delta_\theta^\eps(x,y)$ and diverges for
  $s<\delta_\theta^\eps(x,y)$. If $Q^{s,\eps}_{\theta}(x,y)$ converges for all $s\in\RR$, then we set $\delta_\theta^\eps(x,y)=-\infty$. It is clear that for any $\eps>0$ we have\break $\delta_\theta^\eps(x,y)\le \delta(\Gamma)$, the critical exponent of the Poincar{\'e} series. Unfortunately, since the summation is only over a subset of $\Gamma$, this number may depend on $x$ and $y$. If $\eps>\sqrt2$, then $\Gamma(x,y;\theta,\eps)=\{\gamma\in\Gamma: \gamma y\ne  x\}$; 
 by discreteness of $\Gamma$ we have $\gamma y\ne x$ for only finitely many $\gamma\in\Gamma$, hence in this case $\delta_\theta^\eps(x,y)= \delta(\Gamma)$. 

An easy calculation shows that using for $n\gg 1$ the definition 
$$N_\theta^\eps(x,y;n):=\#\{ \gamma\in\Gamma:\, \gamma y\ne x,\ d(x,\gamma y)<n,\ \Vert\thet(x,\gamma y)-\theta\Vert<\eps\}\,$$
we have 
\begin{equation}\label{defbylimsup}
\delta_{\theta}^{\eps}(x,y)=\limsup_{n\to\infty}\frac{\log
  N_{\theta}^{\eps}(x,y;n)}{n}\,.
  \end{equation}
\begin{df}
The number $\delta_{\theta}(\Gamma):=\displaystyle\lim_{\eps\to
  0}\delta_{\theta}^{\eps}(o,o)$ is called the {\hd exponent of
  growth of $\Gamma$ of slope $\theta$}.
\end{df}
Notice that the exponent of growth  does not depend on the choice of arguments of
$\delta_{\theta}^\eps$ by elementary geometric estimates. 
Moreover, at first sight this definition for $\delta_\theta(\Gamma)$ seems to be different from the one given  in the introduction using $N_\theta^\eps(n)$; however, since for any unit vector $\theta\in\RR_{\ge 0}^r$, all $\eps>0$  and $n\in\NN$ one has
\[ N_\theta^\eps(x,x;n)\le N_\theta^\eps(n)\le N_\theta^{\eps\sqrt{r}}(x,x;n)\]
the definitions obviously are equivalent.

Before we state properties of the exponent of growth, we illustrate the notion by means of an important

\noindent{\sc Example:}$\quad$
Suppose $\XX$ is a product $\XX=\XX_1\times \XX_2\times\cdots\times\XX_r$ of Hadamard manifolds with pinched negative curvature, and 
$\Gamma_i<\is(\XX_i)$  is a convex cocompact
group with critical exponent $\delta_i>0$, $i\in\{1,2,\ldots,r\}$. The product group 
$\Gamma=\Gamma_1\times \Gamma_2\times\cdots\times\Gamma_r$ 
then acts on the product manifold $\XX$, and for any unit vector 
$ \theta\in\RR^r_{\ge 0}$ with coordinates $\theta_1,\theta_2,\ldots,\theta_r\ge 0$ 
we have
\[\delta_\theta(\Gamma)= \sum_{i=1}^r\delta_i \theta_i.\]
\prf\  By 
Theorem~6.2.5 in \cite{MR1348871} (compare also \cite{MR1007619}) there exists a constant $C>0$ \st for all
 $i\in\{1,2,\ldots, r\}$ and $\rho_i>0$ on has 
\begin{eqnarray*}\label{groest}\frac1{C} \e^{\delta_i R} & \le & \#\{\gamma_i\in\Gamma_i : \ R-\rho_i\le d_i(\xo_i,\gamma_i \xo_i)<R\}\nonumber\\
&\le &\#\{\gamma_i\in\Gamma_i :
d_i(\xo_i,\gamma_i \xo_i)<R\}\le C \e^{\delta_i R}\end{eqnarray*}
if $R\gg 1$ is sufficiently large.
 Given $\theta\in E$, we estimate for $\eps>0$ sufficiently small and $n>\sqrt{r}/\eps$ the number of orbit points 
\begin{align*}
\Delta
N_\theta^\eps(\xo,\xo;n)&= \#\{\gamma=(\gamma_1,\gamma_2,\ldots,\gamma_r)\in\Gamma:\ \Vert\thet(\xo,\gamma\xo)-\theta\Vert<\eps, \\
&\hspace*{1.3cm} n-1\le \sqrt{d_1(\xo_1,\gamma_1 \xo_1)^2+d_2(\xo_2,\gamma_2\xo_2)^2+\cdots+d_r(\xo_r,\gamma_r\xo_r)^2} 
<n\}. 
\end{align*} 
We first notice that if $\theta_i(\xo,\gamma\xo)\in [0,1]$, $i\in\{1,2,\ldots,r\}$, denote the coordinates of $\thet(\xo,\gamma\xo)$, then 
\[ d_i(\xo_i,\gamma_i\xo_i)=d(\xo,\gamma\xo)\cdot \theta_i(\xo,\gamma\xo) \quad\text{for all}\quad  i\in\{1,2,\ldots,r\}.\] 
Moreover, 
$\Vert \thet(\xo,\gamma\xo)-\theta\Vert <\eps$ implies $|\theta_i(\xo,\gamma\xo)-\theta_i|<\eps$ 
for all $i\in\{1,2,\ldots,r\}$. So in particular 
\begin{align*}
\Delta
N_\theta^\eps(\xo,\xo;n)&\le \#\big\{\gamma=(\gamma_1,\gamma_2,\ldots,\gamma_r)\in\Gamma:\  |\theta_i(\xo,\gamma\xo)-\theta_i|<\eps\quad\text{and} \\
& 
\hspace*{3.3cm}  d_i(\xo_i,\gamma_i \xo_i)<n\cdot \theta_i(\xo,\gamma\xo)\quad\text{for all}\quad  i\in\{1,2,\ldots,r\} \big\}\\
&\le   \prod_{i=1}^r \#\{\gamma_i\in\Gamma_i:\ d_i(\xo_i,\gamma_i \xo_i)<n(\theta_i+\eps)\}\\
&\le C^r  \e^{n \big( (\theta_1+\eps)\delta_1 + (\theta_2+\eps)\delta_2+\cdots+(\theta_r+\eps)\delta_r\big)}
\end{align*}
for $n$ sufficiently large.

For the lower bound we first denote $I^+\subset\{1,2,\ldots,r\}$ the set of indices $i$  \st $\theta_i>0$, and $I^o=\{1,2,\ldots,r\}\setminus I^+$. Notice that if $\gamma=(\gamma_1,\gamma_2,\ldots,\gamma_r)\in\Gamma$ satisfies 
\[ 
n-1 \le \frac{d_i(\xo_i,\gamma_i\xo_i)}{\theta_i}<n\quad\text{for all}\quad  i\in I^+\]
and $\gamma_i=\id\,$ for $i\in I^o$, then 
$\, n-1\le d(\xo,\gamma\xo)<n\, $ and -- since $\theta_i\in (0,1]$ for all\break $i\in I^+$  -- we get
\[\theta_i -\frac1{n}\le \theta_i \cdot \frac{n-1}{n}<\frac{d_i(\xo_i,\gamma_i\xo_i)}{d(\xo,\gamma\xo)}<\theta_i\cdot \frac{n}{n-1}\le \theta_i+\frac1{n-1}\quad\text{for}\quad i\in I^+.\]
If $n>\sqrt{r}/\eps$ this implies 
\[ |\theta_i(\xo,\gamma\xo)-\theta_i| <\frac1{n}<\frac{\eps}{\sqrt{r}}\quad\text{ for all}\quad i\in I^+; \] 
since $d_i(\xo_i,\gamma_i\xo_i)=0=\theta_i\, $ for all $i\in I^o$ we obtain
 $\Vert \thet(\xo,\gamma\xo)-\theta\Vert<\eps.$
So we conclude that for $n>\sqrt{r}/\eps$ sufficiently large 
\begin{align*}
 \Delta N_\theta^\eps(\xo,\xo;n)&\ge\#\{(\gamma_1,\gamma_2,\ldots,\gamma_r)\in\Gamma : \ n-1\le \frac{d_i(\xo_i,\gamma_i\xo_i)}{\theta_i}<n\quad\text{for}\quad i\in I^+\\
&\hspace*{6.4cm} \text{and}\quad \gamma_i=\id\quad\text{for}\quad i\in I^o \}\\
&=\prod_{i\in I^+} \#\{\gamma_i\in\Gamma_i:\ n\theta_i-\theta_i\le d_i(\xo_i,\gamma_i\xo_i)\le n\theta_i\}\\
&\ge\frac1{C^r}\cdot \e^{n \sum_{i\in I^+} \delta_i \theta_i}=\frac1{C^r}\cdot \e^{n( \delta_1 \theta_1+\delta_2 \theta_2+\cdots+\delta_r\theta_r)}.
\end{align*}
So by definition of $\delta_{\theta}(\Gamma)$ we get
$\ \delta_\theta(\Gamma)=\delta_1\theta_1 +\delta_2\theta_2+\cdots+\delta_r\theta_r.$\qed\\[3mm]
The first easy property of the exponent of growth in the general case is  
\begin{lem}\label{nichtleer}
For $\theta\in E$ we have
\[\begin{array}{ccr}\delta_{\theta}(\Gamma) \ge 0 &\text{if} &\Lim\cap \rand_\theta \ne \emptyset, \\[1mm]
\delta_{\theta}(\Gamma)=-\infty&\text{if} &\Lim\cap \rand_\theta = \emptyset.\ea\]  
In particular, $P_\Gamma=\{\theta\in E: \delta_\theta(\Gamma)\ge 0\}$.
\end{lem}
\prf\  Suppose $\Lim\cap \rand_\theta\ne \emptyset$.
Then by Lemma~\ref{prodtopologytheta} for any $\eps>0$ there
exist infinitely many $\gamma\in\Gamma$ \st $\Vert\thet(\xo,\gamma\xo)-\theta\Vert<\eps$. In particular
$$ \sum_{\gamma\in\Gamma(o,o;\theta,\eps)} 1 =
  Q^{0,\eps}_{\theta}(o,o) \quad\mbox{diverges},$$
hence $\delta_{\theta}^\eps(o,o)\ge 0$. We conclude $\delta_{\theta}(\Gamma)=\displaystyle \lim_{\eps\to 0}\delta_{\theta}^{\eps}(o,o)\ge 0$.

If $\Lim\cap \rand_\theta= \emptyset$, then for some $\eps>0$ sufficiently small the number of elements $\gamma\in\Gamma$ with $\Vert \thet(\xo,\gamma\xo)-\theta\Vert <\eps$ is finite; otherwise there would be an accumulation point in $\rand_\theta$. 
In particular, we have
\[ Q^{\;0,\eps}_{\theta}(o,o) =\sum_{\gamma\in
  \Gamma(o,o;\theta,\eps)} 1  = \# \Gamma(\xo,\xo;\theta,\eps) <\infty,\quad\text{i.e.}\quad \delta_{\theta}^\eps(o,o)\le 0.\]
Moreover, if $s\le 0$ and $d:=\max\{d(\xo,\gamma\xo): \gamma\in \Gamma(o,o;\theta,\eps)\}$, then
\[ Q^{\;s,\eps}_{\theta}(o,o) =\sum_{\gamma\in
  \Gamma(o,o;\theta,\eps)} \e^{|s| d(\xo,\gamma\xo)}\le \sum_{\gamma\in
  \Gamma(o,o;\theta,\eps)} \e^{|s|\cdot d}\le \e^{|s|\cdot d}\cdot  \# \Gamma(\xo,\xo;\theta,\eps)<\infty.\]
We conclude 
  $\delta_{\theta}^\eps(o,o)=-\infty$ and therefore  $\delta_{\theta}(\Gamma)=\displaystyle \lim_{\eps\to 0}\delta_{\theta}^{\eps}(o,o)=-\infty$.\qed\\[3mm]
The following proposition states that the map 
\[ E\to\RR,\quad \theta\mapsto\delta_\theta(\Gamma)\]
 is upper semi-continuous.
\begin{prp}\label{upsemcon}
Let $(\theta^{(j)})\subset E$ be a sequence converging to $\theta\in E$. Then
$$ \limsup_{j\to\infty}\,\delta_{\theta^{(j)}}(\Gamma)\le \delta_{\theta}(\Gamma)\,.$$
\end{prp}
\prf\  Let $\eps_0\in (0,1)$. Then $\theta^{(j)}\to\theta$ implies
$\Vert\theta^{(j)}-\theta\Vert<\eps_0/2$ for $j$ sufficiently large. Let
$\eps\in (0,\eps_0/2)$ and $\gamma\in\Gamma(\xo,\xo;\theta^{(j)},\eps)$. 
Then
$$\Vert\thet(\xo,\gamma\xo)-\theta\Vert<\eps+\eps_0/2<\eps_0\,,$$
hence for $j$ sufficiently large  $\Gamma(\xo,\xo;\theta^{(j)},\eps)\subset \Gamma(\xo,\xo;\theta,\eps_0)$.
This shows 
\[ \delta_{\theta^{(j)}}^{\eps}(o,o)\le\delta_{\theta}^
{\eps_0}(o,o),\quad\text{and therefore}\quad
\delta_{\theta^{(j)}}(\Gamma)=\lim_{\eps\to
  0}\,\delta_{\theta^{(j)}}^{\eps}(o,o)\le\delta_{\theta}^{\eps_0}(o,o).\]
We conclude \vspace{-2mm}
\begin{gather}  \limsup_{j\to\infty}
\delta_{\theta^{(j)}}(\Gamma)\le
\delta_{\theta}^{\eps_0}(o,o)\,,\quad \mbox{hence}
\notag\\
\limsup_{j\to\infty} \delta_{\theta^{(j)}}(\Gamma){=}
\lim_
{\eps_0\to 0}\Big(\!\limsup_{j\to\infty} \delta_{\theta^{(j)}}(\Gamma)\!\Big){\le}
\lim_{\eps_0\to 0}\delta_{\theta}^{\eps_0}(o,o
){=}\delta_{\theta}(\Gamma).\tag*{$\square$
}
\end{gather}

For convenience, we will now consider the homogeneous extension of the exponent of growth to a map $\Psi_\Gamma:\RR_{\ge 0}^r\to \RR$.  
Using a special case of a theorem due to J.-F.~Quint, we will prove that this homogeneous extension $\Psi_\Gamma$ is concave, i.e. for any $\x$, $\y\in\RR_{\ge 0}^r$ and $t\in [0,1]$ one has
\[\Psi_\Gamma(t\x+(1-t)\y)\ge t\Psi_\Gamma(\x)+(1-t)\Psi_\Gamma(\y).\] 
In order to state Quint's Theorem, we recall that in a metric space  the ball of radius $t\ge 0$ centered at  $\x$ is denoted $B(\x,t)$. Moreover, we let  $D$ denote the Dirac measure and $\nu_\Gamma:=\sum_{\gamma\in\Gamma}D_{H(\xo,\gamma\xo)}$ the counting measure on $\RR^r_{\ge 0}$, i.e. the image of the counting measure $\sum_{\gamma\in\Gamma}D_{\gamma}$ of $\Gamma$ by the map $\Gamma\to\RR_{\ge 0}^r$, $\gamma\mapsto H(\xo,\gamma\xo)$.  
\begin{thr}\label{concavity}(\cite{MR1933790}, Theorem 3.2.1)
If there exist $s,t,c>0$ \st for any $\x,\y\in\RR^r$ the inequality
\begin{equation}\label{concavegrowth}
\nu_\Gamma(B(\x+\y,s))\ge c\cdot \nu_\Gamma (B(\x,t))\cdot \nu_\Gamma (B(\y,t))
\end{equation}
 holds, then $\Psi_\Gamma$ is concave.
\end{thr}
So we only have to prove inequality~(\ref{concavegrowth}) which will be done with the help of the generic product Proposition~\ref{genprod}: 
\begin{lem}
There exist $s,t,c>0$ \st for any $\x,\y\in\RR^r$ we have
$$\nu_\Gamma(B(\x+\y,s))\ge c\cdot \nu_\Gamma (B(\x,t))\cdot \nu_\Gamma (B(\y,t))\,.$$
\end{lem}
\prf\  Notice that $\nu_\Gamma (B(\x,t))=\#\{\gamma\in\Gamma:\Vert H(\gamma)-\x\Vert <t\}$. Fix $t>0$, set $s=\kappa+2t$ with $\kappa\ge 0$ from  Proposition~\ref{genprod} (a) and denote $C>0$ the inverse of the cardinality of the set $\Lambda\times\Lambda$ from  Proposition~\ref{genprod} (b). We set 
$$P(\Gamma):=\{(\alpha,\beta)\in\Gamma\times\Gamma :\Vert H(\alpha)-\x\Vert <t\,,\ \Vert H(\beta)-\y\Vert <t\}$$
and will show that for all $\x,\y\in\RR^r$ 
$$\#\{\gamma\in\Gamma : \Vert H(\alpha)-\x-\y\Vert<s\}\ge C\cdot \# P(\Gamma)\,.$$
Let $(\alpha,\beta)\in P(\Gamma)$. Then $\gamma:=\pr(\alpha,\beta)\in\Gamma$ satisfies
\be \Vert H(\gamma)-\x-\y\Vert &\le &\Vert H(\gamma)-H(\alpha)-H(\beta)\Vert +\Vert H(\alpha)-\x\Vert+\Vert H(\beta)-\y\Vert\\
&\le& \kappa+t+t=s\,.\ee
Moreover,  Proposition~\ref{genprod} (b) implies that the number of  different elements in $P(\Gamma)$ which can yield the same element in $\{\gamma\in\Gamma : \Vert H(\gamma)-\x-\y\Vert<s\}$ is bounded by $\#(\Lambda\times\Lambda)$.\qed\\[3mm]
So Theorem~\ref{concavity} gives 
\begin{thr}\label{concave}
The function $\Psi_\Gamma$ is concave.
\end{thr}
This finally allows to complete the proof of the first statement in Theorem~D: 
\begin{thr}\label{PGamma=limitcone}
The set of slopes of limit points of $\,\Gamma$ satisfies $P_\Gamma=\ell_\Gamma\cap E$.
\end{thr}
\prf\ For convenience we denote ${\cal P}_\Gamma$ the set of half-lines in $\RR^r_{\ge 0}$ spanned by all vectors $\theta\in P_\Gamma$. By Lemma~\ref{nichtleer} we have
${\cal P}_\Gamma=\{H\in \RR^r_{\ge 0}: \Psi_\Gamma(H)\ge 0\}$, and concavity of $\Psi_\Gamma$ immediately implies that the cone ${\cal P}_\Gamma$ is convex. Moreover, by Proposition~\ref{propertiesofPGamma} we have $P_\Gamma^{reg}={\cal P}_\Gamma\cap E^+=\ell_\Gamma\cap E^+$, and hence ${\cal P}_\Gamma\cap\RR^r_{> 0}=\ell_\Gamma\cap\RR^r_{> 0}$. Since both  ${\cal P}_\Gamma$ and $\ell_\Gamma$ are closed this gives ${\cal P}_\Gamma=\ell_\Gamma$ and the claim follows from 
$P_\Gamma={\cal P}_\Gamma\cap E$. \qed\\[3mm]
With the notation introduced in the proof of the previous theorem we further remark that Lemma~\ref{nichtleer}
implies 
\[ {\cal P}_\Gamma= \{H\in \RR^r_{\ge 0}: \Psi_\Gamma(H)\ge 0\}=\{H\in \RR^r_{\ge 0}: \Psi_\Gamma(H)>-\infty\};\]
so the fact that ${\cal P}_\Gamma=\ell_\Gamma$ terminates the proof of the first statement in Theorem~E of the introduction. 
In order to show the second statement 
we need to apply the following special case of Lemma~4.1.5 in \cite{MR1933790}:
\begin{lem}\label{Dirichtletexpgrowthpositive}(\cite{MR1933790}, Theorem 3.2.1)
Let $\alpha,\beta\in\Gamma$ be independent regular axial isometries and $\phi_{u,v}:\langle\alpha,\beta\rangle^+\to\RR$ the unique homomorphism of semi-groups sending $\alpha$ to $u$ and $\beta$ to $v$ in the additive group $(\RR,+)$. Then for all $u,v>0$ the Dirichlet series
\[ \sum_{\omega\in\langle \alpha,\beta\rangle^+} \e^{-s \phi_{u,v}(\omega)}\] 
has exponent of convergence $\delta>0$.
\end{lem}
With the help of this lemma we can finally deduce
\begin{thr}\label{expgrowthpositive}
$\Psi_\Gamma$ is strictly positive on the relative interior of $\ell_\Gamma$.
\end{thr}
\prf\ As a first step we will show that $\delta(\Gamma)>0$. Concavity of $\Psi_\Gamma $ then implies that there exists $\theta^*\in 
P_\Gamma$ 
\st $\delta_{\theta^*}(\Gamma)=\delta(\Gamma)>0$. Moreover, concavity and upper-semicontinuity of  $\Psi_\Gamma$ 
imply continuity of $\Psi_\Gamma$  on the closed convex cone
\[ \{ H\in\RR_{\ge 0}^r : \Psi_\Gamma (H)\ge 0\}\]
which  is equal to $\ell_\Gamma$ according to Theorem~\ref{PGamma=limitcone}.
We conclude that $\Psi_\Gamma$ is strictly positive on the relative interior of $\ell_\Gamma$. 


Instead of only proving $\delta(\Gamma)>0$ we next show the stronger statement that for  any $\theta\in \ell_\Gamma\cap E$ and  $\eps>0$ we have 
\[\delta^\eps_\theta(\xo,\xo)>0;\]
notice that this also implies $\delta_\theta(\Gamma)\ge 0$ and hence $\Psi_\Gamma(L)\ge 0$ for all $L\in\ell_\Gamma$. 
According to Proposition~\ref{constrfreesubgrps} there exist independent regular axial isometries $\alpha,\beta\in\Gamma$ with\break 
$\Vert \widehat L(\alpha)-\theta\Vert<\eps$, $\Vert \widehat L(\beta)-\theta\Vert<\eps\,$ 
\st every word  $\omega\in\langle \alpha,\beta\rangle^+$ of length $|\omega|$ in the generators $\alpha,\beta$ satisfies
\[ 
\Vert \widehat H(\xo,\omega\xo)-\theta\Vert<\eps  \quad\text{and}\quad \Vert H(\xo,\omega\xo)-\Phi(\omega)\Vert\le \kappa\cdot |\omega|,\]
where $\Phi :\langle \alpha,\beta\rangle^+\to\RR^r$ is the unique homomorphism of semi-groups sending $\alpha$ to $L(\alpha)$ and $\beta$ to $L(\beta)$.
In particular, using the notation introduced in Lemma~\ref{Dirichtletexpgrowthpositive},  we get from the inverse triangle inequality
\[\big|\Vert H(\xo,\omega\xo)\Vert-\Vert \Phi(\omega)\Vert\big| \le \kappa\cdot |\omega| =\phi_{\kappa,\kappa}(\omega);\]
from $\Vert \Phi(\omega)\Vert \le \phi_{l(\alpha),l(\beta)}(\omega)$ we further obtain
\[ d(\xo,\omega\xo)=\Vert H(\xo,\omega\xo)\Vert\le \phi_{l(\alpha),l(\beta)}(\omega) + \phi_{\kappa,\kappa}(\omega)=\phi_{l(\alpha)+\kappa,l(\beta)+\kappa}(\omega)\]
and therefore
\[ Q^{\;s,\eps}_{\theta}(o,o) =\sum_{\gamma\in
  \Gamma(o,o;\theta,\eps)} \e^{-s d(\xo,\gamma\xo)}\ge \sum_{\omega\in\langle \alpha,\beta\rangle^+} \e^{-s d(\xo,\omega\xo)}\ge  \sum_{\omega\in\langle \alpha,\beta\rangle^+} \e^{-s  \phi_{l(\alpha)+\kappa,l(\beta)+\kappa}(\omega)}.\]
%
Since all constants $l(\alpha)$, $l(\beta)$ and $\kappa$ are positive, Lemma~\ref{Dirichtletexpgrowthpositive} implies that $Q^{s,\eps}_{\theta}(o,o)$ diverges for some $s>0$ and hence $ \delta_\theta^\eps(\xo,\xo)>0$. As a corollary we obtain that
\[ \sum_{\gamma\in\Gamma} \e^{-sd(\xo,\gamma\xo)}\]
has exponent of convergence $\delta(\Gamma)>0$.\qed\\

\bibliography{Bibliographie} 

\vspace{0.9cm}
\noindent Gabriele Link\\
Karlsruher Institut f\"ur Technologie (KIT)\\
Kaiserstr. 12\\[1mm]
D-76128 Karlsruhe\\
e-mail:\ gabriele.link@kit.edu

\end{document}